\documentclass[reqno,12pt]{amsart}
\usepackage{amssymb}
\usepackage{amsmath}
\usepackage{amsfonts}

\theoremstyle{plain}

\newtheorem{corollary}{Corollaire}

\newtheorem{definition}{Définition}
\newtheorem{example}{Exemple}

\newtheorem{lemma}{Lemme}

\newtheorem{proposition}{Proposition}
\newtheorem{remark}{Remarque}

\numberwithin{equation}{section}
\input{tcilatex}

\begin{document}
\title{Opr\'{e}rateur absolument continues et interpolation}
\author{Daher Mohammad }
\email{m.daher@orange.fr}

\begin{abstract}
\ \ 
\end{abstract}

\begin{abstract}
Dans la premi\`{e}re partie de ce travail, on \'{e}tudie les op\'{e}rateurs
absolument continus qui sont d\'{e}finis sur des espaces de fonctions au
sens large.
\end{abstract}

\begin{abstract}
Dans la deuxi\`{e}me partie, on \'{e}tablit quelques r\'{e}sultats
concernant les op\'{e}rateurs absolument continus lorsque les espaces des
fonctions (au sens large) sont des espaces d'interpolation.
\end{abstract}

\maketitle

\begin{abstract}
\textsc{Abstract. }In the first part of this work, we study the absolutely
continuous operators which are defined on fuction spaces with wide sense.

In the second part, we show some results concerning the absoltely continuous
operators when the function spaces (with wide sense) are interpolation
spaces.
\end{abstract}

\bigskip Mots Cl\'{e}s:absolument continu

\begin{center}
\textsc{Introduction.}
\end{center}

Soit $X,Y$ deux espaces de Banach. D\'{e}signons par $\mathcal{L}(X,Y)$ les
op\'{e}rateurs born\'{e}s de $X$ \`{a} valeurs dans $Y$ et $K(X,Y)$ le
sous-espace de $\mathcal{L}(X,Y)$ form\'{e}s des op\'{e}rateurs compacts.

Nous introduisons dans la premi\`{e}re partie de ce travail, les op\'{e}%
rateurs absolument continus dans $\mathcal{L}(X,Y)$ lorsque $X$ est un
espace de fonctions au sens large, cette d\'{e}finition coincide avec celle
de \cite{Ben-Sh}, si $X$ est un espace de fonctions mesurables$.$ Dans la
suite, nous donnons des conditions n\'{e}cessaires pour q'un op\'{e}rateur
dans $K(X,Y)$ soit absolument continu.

Dans la deuxi\`{e}me partie, nous \'{e}tudions les op\'{e}rateurs absolument
continus sur les espaces d'interpolation, comme des espaces de fonctions au
sens large. Finalement, nous montrons que $B_{\theta }^{\ast }=,(B_{0}^{\ast
},B_{1}^{\ast })^{\theta }$ isom\'{e}triquement au sens large, pour tout
couple d'interpolation $(B_{0},B_{1})$ tel que $B_{0}\cap B_{1}$ est dense
dans $B_{0}$ et $B_{1}$ et tout $\theta \in \left] 0,1\right[ .$

\section{Op\'{e}rateurs absolument continus}

Soit $Y$ un espace de Banach complexe, $Y^{\ast }$ son dual. Pour $y\in Y$
et $y^{\ast }\in X^{\ast }$ on note $\left\langle y,y^{\ast }\right\rangle
=y^{\ast }(y).$

\begin{definition}
\label{df}Soient $Y$ un espace de Banach, $(\Omega ,\Sigma ,\mu )$ un espace
mesur\'{e} et $\Delta _{Y}:\Sigma \times Y\rightarrow Y$ une application$.$
On dit que $Y$ est un espace des fonctions au sens large sur $(\Omega
,\Sigma ,\mu ,\Delta _{Y})$ si $\Delta _{Y}$ v\'{e}rifie les conditions
suivantes:
\end{definition}

I) $\Delta _{Y}(A\cap B,f)=\Delta _{Y}(A,\Delta _{Y}(f,B)),$ $\forall A,B\in
\Sigma $ et $\forall f\in Y.$

II) $\Delta _{Y}(A,\alpha f+\beta g)=\alpha \Delta _{Y}(A,f)+\beta \Delta
_{Y}(A,g),\forall A\in \Sigma ,$ $\forall f,g\in Y$ et $\forall \alpha
,\beta \in \mathbb{C}.$

III) $\Delta _{Y}(A\cup B,f)=\Delta _{Y}(A,f)+\Delta _{Y}(B,f),$ $\forall
f\in Y$ et $\forall A,B\in \Sigma $ tel que $A\cap B=\emptyset $ .

IV) Si $\mu (A)=0$ ($A\in \Sigma ),$ alors $\Delta _{Y}(A,f)=0$ et $\Delta
_{Y}(\Omega ,f)=f$ , $\forall f\in Y.$

V) Il existe une constante $C>0$ telle que $\left\Vert \Delta
_{Y}(A,f)\right\Vert \leq C\left\Vert f\right\Vert ,\forall f\in Y$ et $%
\forall A\in \Sigma .$

\begin{example}
\label{cc}Soit $Y$ un espace de fonctions mesurables sur $\left[ 0,\infty %
\right[ $ (cf.\cite{Ben-Sh}). On d\'{e}finit $\Delta _{Y}$ par $\Delta
_{Y}(A,f)=f\mathcal{X}_{A},$ $(A,f)\in \Sigma \times Y.$ Il est facile de
voir que $\Delta _{Y}$ v\'{e}rifie les conditions I),II),III),IV),V).
\end{example}

Pour $(A,f)\in \Sigma \times Y$ on note $\Delta _{Y}(A,f)=f\mathcal{X}_{A}.$

Soit $Y$ un espace de fonctions au sens large sur $(\Omega ,\Sigma ,\mu
,\Delta _{Y})$. Pour $f\in Y$ notons $N(f)=\sup \left\{ \left\Vert f\mathcal{%
X}_{A}\right\Vert ;\text{ }A\in \Sigma \right\} .$ Il est \'{e}vident que $%
N(.)$ est une norme \'{e}quivalente sur $Y$ et $N(f\mathcal{X}_{A})\leq N(f)$
pour tout $(A,f)\in \Sigma \times Y.$ On peut donc supposer que $C=1$ dans
V).

Soit $Y$ un espace de fonctions au sens large sur $(\Omega ,\Sigma ,\mu
,\Delta _{Y})$. \ Pour $(A,f^{\ast })\in \Sigma \times Y^{\ast },$ on d\'{e}%
finit $\Delta _{Y^{\ast }}^{\ast }(A,f^{\ast })\in Y^{\ast }$ par $%
\left\langle f,\Delta _{Y^{\ast }}^{\ast }(A,f^{\ast })\right\rangle
=\left\langle \Delta _{Y}(A,f),f^{\ast }\right\rangle $ pour tout $f\in Y.$

Nous allons la proposition \'{e}vidente suivante:

\begin{proposition}
\label{hh}Soit $Y$ un espace de fonctions au sens large sur $(\Omega
\,,\Sigma ,\mu ,\Delta _{Y}).$ Alors $Y^{\ast }$ est un espace de fonctions
au sens large sur $(\Omega ,\Sigma ,\mu ,\Delta _{Y^{\ast }}^{\ast }).$
\end{proposition}

Soient $Y$ un espace de Banach et $Z$ un sous-espace ferm\'{e} de $X.$ On
note $\left[ x\right] $ l'image de $x$ dans l'espace quotient $X/Y.$

\begin{definition}
\label{bn}Soient $Y$ un espace de fonctions au sens large sur ($\Omega
,\Sigma ,\mu ,\Delta _{Y})$ et $Z$ un sous-espace de Banach de $Y$. On dit
que $Z$ est stable par $\Delta _{Y}$, si $\Delta _{Y}(A,f)\in Z$ pour tout $%
(A,f)\in \Sigma \times Z.$
\end{definition}

Consid\'{e}rons $Y$ un espace de fonctions au sens large sur ($\Omega
,\Sigma ,\mu ,\Delta _{Y})$ et $Z$ un sous-espace stable par $\Delta _{Y}$.
On d\'{e}finit $\Delta _{X/Y}:\Sigma \times Y/Z\rightarrow Y/Z,$ par $\Delta
_{Y/Z}(A,\left[ f\right] )=\left[ f\mathcal{X}_{A}\right] $ pour tout $A\in
\Sigma $ et tout $f\in Y.$

\begin{proposition}
\label{kk}Soient $Y$ est un espace de fonctions au sens large sur ($\Omega
,\Sigma ,\mu ,\Delta _{Y})$ et $Z$ un sous-espace de Banach de $Y$ stable
par $\Delta _{Y}$. Alors $Y/Z$ est un espace de fonctions au sens large sur $%
(\Omega ,\Sigma ,\mu ,\Delta _{Y/Z}).$
\end{proposition}

D\'{e}monstration.

Il est facile de voir que $\Delta _{Y/Z}$ v\'{e}rifie I), II), III),IV).
Montrons que $\Delta _{Y/Z}$ v\'{e}rifie la condition V).

Pour tout $(A,f)\in \Sigma \times Y$ et tout $g\in \left[ f\right] $ on a $%
\left\Vert \Delta _{Y/Z}(A,\left[ f\right] )\right\Vert _{Y/Z}=\left\Vert %
\left[ \Delta (A,f\right] )\right\Vert _{Y/Z}\leq \left\Vert \Delta
(A,g)\right\Vert _{Y}\leq \left\Vert g\right\Vert .$ Par cons\'{e}quent $%
\left\Vert \Delta _{Y/Z}(A,\left[ f\right] )\right\Vert _{Y/Z}\leq
\left\Vert \left[ f\right] \right\Vert _{Y/Z}$.$\blacksquare $

Soient $X,Y$ deux espaces de fonctions au sens large sur $(\Omega ,\Sigma
,\mu ,\Delta _{X})$,$(\Omega ,\Sigma ,,\mu ,\Delta _{Y})$ respectivement et $%
T:X\rightarrow Y$ est un op\'{e}rateur born\'{e}. \ L'identification $%
T(\Delta _{X})=\Delta _{Y}$ signifie que $T\left[ \Delta _{X}(A,f)\right]
=\Delta _{Y}(A,T(f))$ pour tout $(A,f)\in \Sigma \times X.$

\begin{definition}
\label{hm}Soit $Y$ un espace de fonctions au sens large sur $(\Omega ,\Sigma
,\mu ,\Delta _{Y}):$
\end{definition}

a) \ Soit $f\in Y$. On dit que $f$ est absolument continu, si pour tout $%
\varepsilon >0,$ il existe $\delta >0$ tel que si $\mu (A)<\delta ,$ alors $%
\left\Vert f\mathcal{X}_{A}\right\Vert _{Y}<\varepsilon .$

b) Une partie $Z$ de $Y$ est dite absolument continue, si pour tout $f\in Z$%
, $\ f$ est absolument continu.

c) Une partie $Z$ de $Y$ est dite uniform\'{e}ment absolument continue, si
pour tout $\varepsilon >0,$ il existe $\delta >0,$ tel que si $\mu
(A)<\delta ,$ alors $\left\Vert f\mathcal{X}_{A}\right\Vert <\varepsilon $
pour tout $f\in Z.$

Soient $Y$ un espace de fonctions au sens large sur ($\Omega ,\Sigma ,\mu
,\Delta _{Y})$ et $Z$ une partie de $Y.$ Posons $V=\ell ^{\infty }(Z,Y).$ On
d\'{e}finit $\Delta _{V}(A,(f_{i})_{i\in Z})=(f_{i}\mathcal{X}_{A})_{i\in
Z}, $ $A\in \Sigma $, $(f_{i})_{i\in Z}\in V.$ Il est clair que $V$ est un
espace de fonctions au sens large sur $(\Omega ,\Sigma ,\mu ,\Delta _{V}).$
Supposons que $Z$ est uniform\'{e}ment born\'{e}. Il est facile de voir que $%
Z$ est uniform\'{e}ment absolument continue, si et seulement si l'\'{e}l\'{e}%
ment $(f)_{f\in Z}$ est absolument continu dans $V.$

\begin{remark}
\label{lr}Soient $Y$ un espace de fonctions au sens large sur $(\Omega
,\Sigma ,\mu ,\Delta _{Y})$ et $Z$ une partie absolument continue dans $Y.$
Alors l'adh\'{e}rence de $Z$ dans $Y$ est absolument continue.
\end{remark}

Pour tout $f\in Y$ notons $\nu _{f}(A)=(f\mathcal{X}_{A}),$ $A\in \Sigma ,$
la mesure $\nu _{f}$ est finiment additive \`{a} valeurs dans $Y.$

\begin{definition}
\label{hy}Soient $Y$ un espace de fonctions au sens large sur $(\Omega
,\Sigma ,\mu ,\Delta _{Y})$ et $f\in Y.$ On dit que la mesure $\nu _{f}$ est 
$\mu -$d\'{e}nombrablement additive, si pour toute suite $(A_{k})_{k\geq 0}$
dans $\Sigma ,$ deux-\`{a}-deux disjoints v\'{e}rifiant $\mu (\underset{%
k\geq 0}{\cup }A_{k})<+\infty ,$ alors $\nu _{f}(\underset{k\geq 0}{\cup }%
A_{k})=\underset{k\geq 0}{\dsum }\nu _{f}(A_{k}).$
\end{definition}

\begin{lemma}
\label{ez}Soient $Y$ un espace de fonctions au sens large sur $(\Omega
,\Sigma ,\mu ,\Delta _{Y})$ et $f\in Y.$ Les assertions suivantes sont \'{e}%
quivalentes.
\end{lemma}

1) $f$ est absolument continu.

2) \ $\nu _{f}$ est $\mu -$d\'{e}nombrablement additive.

3) \ Pour toute suite d\'{e}croissante $(A_{n})_{n\geq 0}$ dans $\Sigma $
telle que $\mu (A_{n})\underset{n\rightarrow \infty }{\rightarrow }0,$ alors 
$\nu _{f}(A_{n})\underset{n\rightarrow \infty }{\rightarrow }0$ dans $Y.$

D\'{e}monstration.

1) $\Longrightarrow $ 2).

Fixons $\varepsilon >0.$ Il existe $\delta >0$ tel que si $\mu (A)<\delta ,$ 
$\left\Vert (f\mathcal{X}_{A})\right\Vert _{Y}<\varepsilon .$

Soient $(A_{k})_{k\geq 0}$ une suite dans $\Sigma $ deux-\`{a}-deux
disjoints telle que $\mu (\underset{k\geq 0}{\cup }A_{k})<\infty .$ Il
existe $n_{0}\geq 1$ tel que $\mu (\underset{k\geq n}{\cup }A_{k})=\underset{%
k\geq n}{\dsum }\mu (A_{k})<\delta $ pour tout $n\geq n_{0}.$ Donc

\bigskip

\begin{eqnarray*}
&&\left\Vert \nu _{f}(\underset{k\geq 0}{\cup }A_{k})-\underset{k\leq n}{%
\dsum }\nu _{f}(A_{k})\right\Vert \\
&=&\left\Vert \nu _{f}(\underset{k\geq 0}{\cup }A_{k})-\nu _{f}(\underset{%
k\leq n}{\cup }A_{k})\right\Vert \\
&=&\left\Vert \nu _{f}(\underset{k>n}{\cup A_{k})}\right\Vert \\
&=&\left\Vert (f\mathcal{X}\underset{k>n}{_{\cup }A_{k}})\right\Vert
<\varepsilon .
\end{eqnarray*}%
Par cons\'{e}quent $\nu _{f}(\underset{k\geq 0}{\cup }A_{k}\underset{}{)%
\overset{}{=\underset{k\geq 0}{\dsum }}}\nu _{f}(A_{k}).\blacksquare $

2) $\Longrightarrow $ 3).

Soit $(A_{n})_{n\geq 0}$ une suite d\'{e}croissante dans $\Sigma $ telle que 
$\mu (A_{n})\underset{n\rightarrow +\infty }{\rightarrow }0.$ Il existe une
suite $(A_{n_{k}})_{k\geq 0}$ telle que $\mu (A_{n_{k}})<2^{-k},$ $\forall
k\in \mathbb{N}.$ Notons pour tout $k\in \mathbb{N}$ $%
B_{k}=A_{n_{k}}-A_{n_{k+1}}.$ Il est clair que $\mu (\underset{k\geq 0}{\cup 
}B_{k})<+\infty ,$ donc $\nu _{f}(\underset{k\geq 0}{\cup }B_{k})=\underset{%
k\geq 0}{\dsum }\nu _{f}(B_{k}).$ Comme $\mu (\underset{j\geq k+1}{\cap }%
A_{n_{j}})=0$ d'apr\`{e}s la condition IV), $\nu _{f}(\underset{j\geq k+1}{%
\cap }A_{n_{j}})=0$ pour tout $k.$ D'autre part, $A_{n_{k}}=\left[ \underset{%
j\geq k}{\cup }B_{j}\right] \cup \left[ \underset{j\geq k+1}{\cap }A_{n_{j}}%
\right] ,$ ceci implique que $\nu _{f}(A_{n_{k}})=\nu _{f}(\underset{j\geq k}%
{\cup }B_{j})=\underset{j\geq k}{\dsum }\nu _{f}(B_{j})$ pour tout $k.$
Fixons $\varepsilon >0,$ il existe $k_{0}\in \mathbb{N}$ v\'{e}rifiant $%
\left\Vert \nu _{f}(A_{n_{k_{0}}})\right\Vert =\left\Vert \nu _{f}(\underset{%
k\geq k_{0}}{\cup }B_{k})\right\Vert <\varepsilon .$ Il en r\'{e}sulte que $%
\left\Vert \nu _{f}(A_{n})\right\Vert \leq \left\Vert \nu
_{f}(A_{k_{k_{0}}})\right\Vert <\varepsilon $ pour tout $n\geq
n_{k_{0}}.\blacksquare $

3) $\Longrightarrow $ 1).

Supposons qu'il existe $\varepsilon _{0}$ tel que pour tout $k\in \mathbb{N}%
, $ il existe $B_{k}\in \Sigma $, v\'{e}rifiant $\mu (B_{k})<2^{-k}$ et $%
\left\Vert \nu _{f}(B_{k})\right\Vert >\varepsilon _{0}.$ Observons que $\mu
(\underset{k\geq n}{\cup }B_{k})\leq \underset{k\geq n}{\dsum }\mu
(B_{k})<+\infty .$ Notons pour tout $n\in \mathbb{N}$ $A_{n}=\underset{k\geq
n}{\cup }B_{k}$. Nous avons alors $\mu (A_{n})\underset{n\rightarrow +\infty 
}{\rightarrow }0,$ donc il existe $n_{0}$ tel que $\left\Vert \nu
_{f}(A_{n_{0}})\right\Vert <\varepsilon _{0}.$

D'autre par, $B_{n_{0}}\subset A_{n_{0}},$ par cons\'{e}quent $\left\Vert
\nu _{f}(B_{n_{0}})\right\Vert \leq \left\Vert \nu
_{f}(A_{n_{0}})\right\Vert <\varepsilon _{0},$ ce qui est impossible. Il en r%
\'{e}sulte que $f$ est absolument continu.$\blacksquare $

\begin{remark}
\label{dd}Dans le lemme \ref{ez}, on a toujours 1) $\Longrightarrow $ 2)
sans la condition $V)$ de la d\'{e}finition \ref{df}.
\end{remark}

Soient $X$ un espace de fonctions au sens large sur $(\Omega ,\Sigma ,\mu
,\Delta _{X})$ et $Y$ un espace de Banach. D\'{e}finissons $\Delta _{%
\mathcal{L}(X,Y)}(A,T)(f)=T\mathcal{X}_{A}(f)=T(f\mathcal{X}_{A}),$ $T\in 
\mathcal{L}(X,Y),A\in \Sigma ,f\in X.$ Il est clair que $\Delta _{\mathcal{L}%
(Y,X)}$ v\'{e}rifie les conditions $($I),(II),(III),IV$)$,V), par cons\'{e}%
quent $\mathcal{L}(Y,X)$ est un espace de fonctions au sens large sur $%
(\Omega ,\Sigma ,\mu ,\Delta _{\mathcal{L}(X,Y)}).$

Soient $X$ un espace de Banach et $Y$ un espace de fonctions au sens large
sur $(\Omega ,\Sigma ,\mu ,\Delta _{Y})$. D\'{e}finissons $\Delta _{\mathcal{%
L}^{\prime }(X,Y)}^{\prime }(A,T)(f)=T\mathcal{X}_{A}(f)=(T(f))\mathcal{X}%
_{A}=\Delta _{Y}(A,T(f)),$ $T\in \mathcal{L}(X,Y),A\in \Sigma ,f\in X.$ On v%
\'{e}rifie facilement que $\Delta _{\mathcal{L}(Y,X)}^{\prime }$ v\'{e}rifie
les conditions $($I),(II),(III),IV$)$,V), donc $\mathcal{L}(Y,X)$ est un
espace de fonctions au sens large sur $(\Omega ,\Sigma ,\mu ,\Delta _{%
\mathcal{L}(X,Y)}^{\prime }).$

\begin{definition}
\label{dg}Soient $X$ un espace de fonctions au sens large sur $(\Omega
,\Sigma ,\mu ,\Delta _{X}),$ $Y$ un espace de Banach et $T:X\rightarrow Y$
un op\'{e}rateur born\'{e}. On dit que $T$ est absolument continu, si $T$
est absolument continu par rapport \`{a} $(\Omega ,\Sigma ,\mu ,\Delta _{%
\mathcal{L}(X,Y)})$.
\end{definition}

\begin{definition}
\label{pm}Soient $X$ un espace de Banach, $Y$ un espace de fonctions au sens
large sur $(\Omega ,\Sigma ,\mu ,\Delta _{Y})$ et $T:X\rightarrow Y$ un op%
\'{e}rateur born\'{e}. On dit que $T$ est absolument continu \`{a} gauche,
si $T$ est absolument continu par rapport \`{a} $(\Omega ,\Sigma ,\mu
,\Delta _{\mathcal{L}(X,Y)}^{\prime }).$
\end{definition}

Il est facile de montrer le lemme suivant:

\begin{lemma}
\bigskip \label{rh}Soient $X$ un espace de fonctions au sens large sur ($%
\Omega ,\Sigma ,\mu ,\Delta _{X})$ et $T:X\rightarrow Y$ un op\'{e}rateur
born\'{e}. Alors $T$ est absolument continu si et seulment si $T^{\ast
}:Y^{\ast }\rightarrow X^{\ast }$ est absolument continu \`{a} gauche.
\end{lemma}

\begin{remark}
\label{nh}Soient $X,Y$ deux espaces de fonctions au sens large sur ($\Omega
,\Sigma ,\mu ,\Delta _{X}),(\Omega ,\Sigma ,\mu ,\Delta _{Y})$
respectivement et $T:X\rightarrow Y$ un op\'{e}rateur born\'{e}. Si $%
T(\Delta _{X})=\Delta _{Y},$ alors $T$ est absolument continu si et
seulement si $T$ est absolument continu \`{a} gauche.
\end{remark}

\begin{remark}
\bigskip \label{je}Soient $Y$ un espace de fonctions au sens large sur $%
(\Omega ,\Sigma ,\mu ,\Delta _{Y})$ et $K$ un compact dans $Y$ pour la norme$%
.$ Si $K$ est absolument continue, alors $K$ est uniform\'{e}ment absolument
continue.
\end{remark}

Preuve. En effet,

Consid\'{e}rons $(B_{n})_{n\geq 0}$ une suite d\'{e}croissante dans $\Sigma $
telle que $\mu (B_{n})\underset{n\rightarrow \infty }{\rightarrow }0.$ Soit $%
\ f\in K.$ Il existe $m=m(f)$ tel que $\left\Vert f\mathcal{X}%
_{B_{m}}\right\Vert _{Y}<\varepsilon .$ Notons pour tout $m\in \mathbb{N}$ $%
O_{m}=\left\{ f\in K;\text{ }\left\Vert f\mathcal{X}_{B_{m}}\right\Vert
_{Y}<\varepsilon \right\} $. D'apr\`{e}s ce qui pr\'{e}c\`{e}de $%
(O_{m})_{m\geq 0}$ est un recouvrement ouvert de $K,$ comme $K$ est un
compact, il existe $m_{1},...,m_{m}\in \mathbb{N}$ verifiant $K\subset
O_{m_{1}}\cup ...\cup O_{m_{m}}.$ Noter $m_{0}=\max (m_{1},...,m_{m})$.
Montrons que $\left\Vert f\mathcal{X}_{B_{n}}\right\Vert _{Y}<\varepsilon $
pour tout $f\in K$ et pour tout $n\geq m_{0}.$

Soient $n\geq m_{0}$ et $f\in K.$ Il existe $m_{j}$ tel que $\left\Vert \nu
_{f}(B_{m_{j}})\right\Vert _{Y}<\varepsilon .$ On a alors $\left\Vert \nu
_{f}(B_{n})\right\Vert _{Y}\leq \left\Vert \nu _{f}(B_{m_{j}})\right\Vert
_{Y}<\varepsilon ,$ car $B_{n}\subset B_{m_{j}}.$ Donc sup$_{f\in
K}\left\Vert \nu _{f}(B_{n})\right\Vert _{Y}<\varepsilon .$ D'apr\`{e}s le
lemme \ref{ez}, $(f)_{f\in K}$ est absolument continu dans $\ell ^{\infty
}(K,X)$. Finalement d'apr\`{e}s la remarque \ref{nh}, $K$ est uniform\'{e}%
ment absolument continue.$\blacksquare $

$\blacksquare $

\begin{definition}
\label{bv}Soient $X$ un espace de fonctions au sens large sur $(\Omega
,\Sigma ,\mu ,\Delta _{X})$ et $T:X\rightarrow Y$ un op\'{e}rateur born\'{e}%
. On dit que $T$ est ponctuellement faiblement absolument continu (par
rapport \`{a} ($\Omega ,\Sigma ,\mu ,\Delta _{X}))$, si pour tout $f\in X$,
tout $y^{\ast }\in Y^{\ast }$ et tout $\varepsilon >0$, il existe $\delta >0$
tels que si $\mu (A)<\delta ,$ alors $\left\vert \left\langle T(f\mathcal{X}%
_{A}),y^{\ast }\right\rangle \right\vert <\varepsilon .$
\end{definition}

Soient $X$ un espace de fonctions au sens large sur $(\Omega ,\Sigma ,\mu
,\Delta _{X})$ et $T:X\rightarrow Y$ un op\'{e}rateur born\'{e}. Pour tout $%
f\in X$ on d\'{e}finit la mesure $\nu _{T(f)}$ \ par $\nu _{T(f)}(A)=T(f%
\mathcal{X}_{A}).$

\begin{lemma}
\label{tg}Soient $X$ un espace de fonctions au sens large sur $(\Omega
,\Sigma ,\mu ,\Delta _{X})$ et $T:X\rightarrow Y$ un op\'{e}rateur born\'{e}%
. Supposons que $T$ soit ponctuellement faiblement absolument continu. Alors
pour tout $f\in X$ $\nu _{T(f)}$ est $\mu -$d\'{e}nombrablement additive$.$
\end{lemma}

D\'{e}monstration.

Soit $(A_{k})_{k\geq 0}$ une suite de sous-ensembles mesurables deux-\`{a}%
-deux disjoints telle que $\underset{k\geq 0}{\dsum }\mu (A_{k})<+\infty $.
Comme $T$ est ponctuellement faiblement absolument continu, d'apr\`{e}s la
remarque \ref{dd} pour tout $M\subset \mathbb{N},$ la serie $\underset{k\in M%
}{\dsum }T(f\mathcal{X}_{A_{k}})$ converge faiblement dans $Y$ vers $T(f%
\mathcal{X}\underset{k\geq 0}{_{\cup A_{k}}}).$ D'apr\`{e}s \cite[%
Coroll.4,Chap.I-4]{DU}, la s\'{e}rie$\underset{k\geq 0}{\dsum }T(f\mathcal{X}%
_{A_{k}})$ converge inconditionnellement dans $Y$ vers $T(f\mathcal{X}_{%
\underset{k\geq 0}{\cup }A_{k}}).\blacksquare $

\begin{definition}
\label{vx}Soient $X$ un espace de fonctions au sens large sur $(\Omega
,\Sigma ,\mu ,\Delta _{X}),$ $Y$ un espace de Banach et $T:X\rightarrow Y$
un op\'{e}rateur born\'{e}. On dit que $T^{\ast \ast }$ est ponctuellement pr%
\'{e}faiblement absolument continu, si la mesure $A\in \Sigma \rightarrow
\left\langle y^{\ast },T^{\ast \ast }(f^{\ast \ast }\mathcal{X}%
_{A})\right\rangle $ =$\left\langle y^{\ast },T^{\ast \ast }(\Delta
_{X^{\ast \ast }}^{\ast \ast }(A,f^{\ast \ast }))\right\rangle =\left\langle
y^{\ast },\nu _{T^{\ast \ast }}(A)(f^{\ast \ast })\right\rangle $ est $\mu -$%
d\'{e}nombrablement additive, pour tout $f^{\ast \ast }\in X^{\ast \ast }$.
\end{definition}

Pour tout Banach $X$ on note $B_{X}$ la boule unit\'{e} ferm\'{e}e de $X.$

\begin{lemma}
\label{kg}Soient $X$ un espace de fonctions au sens large sur $(\Omega
,\Sigma ,\mu ,\Delta _{X}),$ $Y$ un espace de Banach et $T:X\rightarrow Y$
un op\'{e}rateur born\'{e}. Alors $\left\langle y^{\ast },T^{\ast \ast
}(\Delta _{X^{\ast \ast }}^{\ast \ast }(A,f^{\ast \ast }))\right\rangle
=\left\langle (T\mathcal{X}_{A})^{\ast }y^{\ast },f^{\ast \ast
}\right\rangle $ pour tout $(A,y^{\ast })\in \Sigma \times Y^{\ast }$ et
tout $f^{\ast \ast }\in X^{\ast \ast }.$
\end{lemma}

D\'{e}monstration.

Soit $(A,y^{\ast })\in \Sigma \times Y^{\ast }.$ Il est clair que
l'application $f^{\ast \ast }\rightarrow \left\langle y^{\ast },T^{\ast \ast
}(\Delta _{X^{\ast \ast }}^{\ast \ast }(A,f^{\ast \ast }))\right\rangle
=\left\langle T^{\ast }y^{\ast },\Delta _{X^{\ast \ast }}^{\ast \ast
}(A,,f^{\ast \ast })\right\rangle =\left\langle \Delta _{X^{\ast }}^{\ast
}(A,T^{\ast }(y^{\ast })),f^{\ast \ast })\right\rangle $ est pr\'{e}%
faiblement continue : $B_{X^{\ast \ast }}\rightarrow \mathbb{C}$ et
l'application $f^{\ast \ast }\rightarrow \left\langle (T\mathcal{X}%
_{A})^{\ast }y^{\ast },f^{\ast \ast }\right\rangle $ est pr\'{e}faiblement
continue:$B_{X^{\ast \ast }}\rightarrow \mathbb{C}$. D'autre part, $%
\left\langle y^{\ast },T^{\ast \ast }(\Delta _{X^{\ast \ast }}^{\ast \ast
}(A,f))\right\rangle =\left\langle y^{\ast },T^{{}}(\Delta
_{X}(A,f))\right\rangle $=$\left\langle (T\mathcal{X}_{A})^{\ast }y^{\ast
},f\right\rangle $, pour tout $f\in B_{X}.$ Comme $B_{X}$ est pr\'{e}%
faiblement dense dans $B_{X^{\ast \ast }},$ alors $\left\langle y^{\ast
},T^{\ast \ast }(\Delta _{X^{\ast \ast }}^{\ast \ast }(A,f^{\ast \ast
}))\right\rangle =\left\langle (T\mathcal{X}_{A})^{\ast }y^{\ast },f^{\ast
\ast }\right\rangle ,$ pour tout $f^{\ast \ast }\in X^{\ast \ast
}.\blacksquare $

\begin{lemma}
\label{ssd}Soient $X$ un espace de fonctions au sens large sur $(\Omega
,\Sigma ,\mu ,\Delta _{X})$ $,$ $Y$ un espace de Banach et $T:X\rightarrow Y$
un op\'{e}rateur ponctuellement faiblement absolument continu. Supposons que 
$\mu $ est une mesure born\'{e}e, $X$ est s\'{e}parable et ne contient pas $%
\ell ^{1}$ isomorphiquement. Alors $T^{\ast \ast }$ est ponctuellement pr%
\'{e}faiblement absolument continu.
\end{lemma}

D\'{e}monstration.

Soient $f^{\ast \ast }\in X^{\ast \ast }$ et $y^{\ast }\in Y^{\ast }$. Comme 
$X$ ne contient pas $\ell ^{1}$ isomorphiquement, d'apr\`{e}s \cite{Ros}, il
existe une suite born\'{e}e $(f_{n})_{n\geq 0}$ dnas $X$ telle que $f_{n}%
\underset{n\rightarrow \infty }{\rightarrow }f^{\ast \ast }$ pr\'{e}%
faiblement dans $X^{\ast \ast }.$ Pour tout $n$ on d\'{e}finit la mesure $%
\nu _{n}$ par $\nu _{n}(A)=\left\langle T(f_{n}\mathcal{X}_{A}),y^{\ast
}\right\rangle ,$ $A\in \Sigma .$ Remarquons d'apr\`{e}s le lemme \ref{kg}
que $\nu _{n}(A)\underset{n\rightarrow \infty }{\rightarrow }\left\langle
y^{\ast },T^{\ast \ast }(f^{\ast \ast }\mathcal{X}_{A})\right\rangle $ pour
tout $A\in \Sigma ,$ donc d'apr\`{e}s \cite[Cor.6,Chap.1-5]{DU}, $\lim_{\mu
(A)\rightarrow 0}\nu _{n}(A)$ existe uniform\'{e}ment en n. Par cons\'{e}%
quent $\lim_{\mu (A)\rightarrow 0}\left\langle y^{\ast },T^{\ast \ast
}(f^{\ast \ast }\mathcal{X}_{A})\right\rangle =0.\blacksquare $

\begin{lemma}
\label{gr}Soient $X$ un espace de fonctions au sens large sur $(\Omega
,\Sigma ,\mu ,\Delta _{X}),$ $Y$ un espace de Banach et $T:X\rightarrow Y$
un op\'{e}rateur ponctuellement faiblement absolument continu. Supposons que 
$X$ est un espace de Grothendieck. Alors $T^{\ast \ast }$ est ponctuellement
pr\'{e}fiablement absolument continu.
\end{lemma}

D\'{e}monstration.

Soient $y^{\ast }\in Y^{\ast }$ et $(A_{k})_{k\geq 0}$ une suite dans $%
\Sigma $ telle que $\mu (A_{k})\underset{n\rightarrow \infty }{\rightarrow }%
0.$ L'op\'{e}rateur $T$ est ponctuellement faiblement continu, donc $T^{\ast
}(y^{\ast })\mathcal{X}_{A_{n}}\underset{n\rightarrow \infty }{\rightarrow }%
0 $ pr\'{e}faiblement dans $X^{\ast },$ comme $X$ est un espace de
Grothendieck, $T^{\ast }(y^{\ast })\mathcal{X}_{A_{n}}\underset{n\rightarrow
\infty }{\rightarrow }0$ faiblement \cite[p.179]{DU} dans $X^{\ast },$ d'apr%
\`{e}s la remarque \ref{dd}, $T^{\ast \ast }$ est ponctuellement pr\'{e}%
fiablement absolument continu.$\blacksquare $

\begin{proposition}
\label{uv}Soient $X$ un espace de fonctions au sens large sur $(\Omega
,\Sigma ,\mu ,\Delta _{X}),$ $Y$ un espace de Banach et $T:X\rightarrow Y$
un op\'{e}rateur compact. Supposons que $T^{\ast \ast }$ est ponctuellement
pr\'{e}faiblement absolument continu. Alors $T$ est absolument continu.
\end{proposition}

D\'{e}monstration.

\bigskip Il suffit de montrer que $\nu _{T}$ est $\mu -$d\'{e}nombrablement
additive, d'apr\`{e}s le lemme \ref{ez}.

Consid\'{e}rons $(A_{k})_{k\geq 0}$ une suite dans $\Sigma ,$ deux-\`{a}%
-deux disjoints. $T^{\ast \ast }$ est ponctuellement pr\'{e}faiblement
absolument continu, donc $\underset{k\in M}{\dsum }\left\langle y^{\ast
},\nu _{T^{\ast \ast }}(A_{k})(f^{\ast \ast })\right\rangle =\left\langle
y^{\ast },\nu _{T^{\ast \ast }}\underset{k\in M}{(\cup }A_{k})(f^{\ast \ast
})\right\rangle ,$ pour tout $M\subset \mathbb{N}.$

D'autre part, d'apr\`{e}s le lemme \ref{kg}, $\underset{k\in M}{\dsum }%
\left\langle y^{\ast },\nu _{T^{\ast \ast }}(A_{k})(f^{\ast \ast
})\right\rangle $ $\underset{k\in M}{=\dsum }\left\langle y^{\ast },(T%
\mathcal{X}_{A_{k}})^{\ast \ast }(f^{\ast \ast })\right\rangle $ et $%
\left\langle y^{\ast },\nu _{T^{\ast \ast }}\underset{k\in M}{(\cup }%
A_{k})(f^{\ast \ast })\right\rangle =\left\langle y^{\ast },(T\mathcal{X}%
\underset{k\in M}{_{\cup A_{k}}})^{\ast \ast }(f^{\ast \ast })\right\rangle .
$ Il en r\'{e}sulte que $\underset{k\in M}{\dsum }\left\langle y^{\ast },(T%
\mathcal{X}_{A_{k}})^{\ast \ast }(f^{\ast \ast })\right\rangle =\left\langle
y^{\ast },(T\mathcal{X}\underset{k\in M}{_{\cup A_{k}}})^{\ast \ast
}(f^{\ast \ast })\right\rangle .$ Comme pour tout $A\in \Sigma ,$ $\nu
_{T}(A)=T\mathcal{X}_{A}$ est un op\'{e}rateur compact $,$ alors d'apr\`{e}s 
\cite[Coroll.3]{Kalt}, $\underset{k\in M}{\dsum }\left\langle \nu
_{T}(A_{k}),u^{\ast }\right\rangle =\left\langle \nu _{T}\underset{k\in M}{%
(\cup }A_{k}),u^{\ast }\right\rangle $ pour tout $M\subset \mathbb{N}$ et
tout $u^{\ast }\in \left[ K(X,Y)\right] ^{\ast }.$ Ceci entra\^{\i}ne d'apr%
\`{e}s \cite[Chap.I-4,Coroll..4]{DU} que la s\'{e}rie $\underset{k\geq 0}{%
\dsum }\nu _{T}(A_{k})$ converge inconditionnellement vers $\nu _{T}(%
\underset{k\geq 0}{\cup }A_{k})$ dans $\mathcal{L}$(X,Y)$.\blacksquare $

\begin{corollary}
\label{gu}Soient $X$, $Y$ deux espaces de fonctions au sens large sur $%
(\Omega ,\Sigma ,\mu ,\Delta _{X}),$ $(\Omega ,\Sigma ,\mu ,\Delta _{Y})$
respectivement tel que $T(\Delta _{X})=\Delta _{Y}$ et $T:X\rightarrow Y$ un
op\'{e}rateur compact. Supposons que $T$ soit ponctuellement faiblement
absolument continu. Alors $T$ est absolument continu.
\end{corollary}

D\'{e}monstration.

D'apr\`{e}s la proposition \ref{uv}, il suffit de montrer que $T^{\ast \ast
} $ est ponctuellement pr\'{e}faiblement absolument continu.

\emph{Etape 1:} Soient $y^{\ast }\in Y^{\ast }$ et $(A_{n})_{n\geq 0}$ une
suite de sous-ensembles mesurables telle que $\mu (A_{n})\underset{%
n\rightarrow \infty }{\rightarrow }0.$ Montrons que $T^{\ast }\left[ \Delta
_{Y^{\ast }}^{\ast }(A_{k},y^{\ast })\right] \rightarrow 0.$

En ffet,

pour tout $f\in X$ on a $\left\langle f,T^{\ast }\left[ \Delta _{Y^{\ast
}}^{\ast }(A_{k},y^{\ast })\right] \right\rangle =\left\langle T(f),\Delta
_{Y^{\ast }}^{\ast }(A_{k},y^{\ast })\right\rangle =\left\langle \Delta
_{Y}(A_{k},T(f)),y^{\ast }\right\rangle =\left\langle T(\Delta
_{X}(A_{k},f)),y^{\ast }\right\rangle .$

Comme $T$ est ponctuellement faiblement absolument continu, $\left\langle
T(\Delta _{X}(A_{k},f)),y^{\ast }\right\rangle \underset{k\rightarrow \infty 
}{\rightarrow }0,$ c'est-\`{a}dire que $T^{\ast }\left[ \Delta _{Y^{\ast
}}^{\ast }(A_{k},y^{\ast })\right] \underset{k\rightarrow \infty }{%
\rightarrow }0$ pr\'{e}faiblement.

D'autre part, $T^{\ast }$ est un op\'{e}rateur compact et la suite $(\Delta
_{Y^{\ast }}^{\ast }(A_{k},y^{\ast }))_{k\geq 0}$ est born\'{e}e, donc $%
T^{\ast }\left[ \Delta _{Y^{\ast }}^{\ast }(A_{k},y^{\ast })\right] \underset%
{k\rightarrow \infty }{\rightarrow }0$ en norme dans $X^{\ast }.$

\emph{Etape 2: }Montrons que $T^{\ast \ast }$ est ponctuellement pr\'{e}%
faiblement absolument continu.

Pour cela, soit $f^{\ast \ast }\in X^{\ast \ast }.$ Remarquons que pour tout 
$y^{\ast }\in Y^{\ast }$ on a%
\begin{eqnarray*}
\left\langle y^{\ast },T^{\ast \ast }(\Delta _{X^{\ast \ast }}^{\ast \ast
}(A_{k},f^{\ast \ast }))\right\rangle &=&\left\langle T^{\ast }(y^{\ast
}),\Delta _{X^{\ast \ast }}^{\ast \ast }(A_{k},f^{\ast \ast })\right\rangle
\\
&=&\left\langle \Delta _{X^{\ast }}^{\ast }(A_{k},T^{\ast }(y^{\ast
})),f^{\ast \ast }\right\rangle =\left\langle T^{\ast }\left[ \Delta
_{Y^{\ast }}^{\ast }(A_{k},y^{\ast })\right] ,f^{\ast \ast }\right\rangle 
\underset{}{\underset{k\rightarrow \infty }{\rightarrow }}0,
\end{eqnarray*}

car $T^{\ast }(\Delta _{Y^{\ast }}^{\ast })=\Delta _{X^{\ast }}^{\ast }$.
Donc $T^{\ast \ast }$ est ponctuellement pr\'{e}faiblement absolument
continu.$\blacksquare $

\begin{remark}
\label{vo}Soit $X$ un espace de Banach au sens large sur $(\Omega ,\Sigma
,\mu ,\Delta _{X}),$ $Y$ un espace de Banach, $X^{\ast }$ est absolument
continu et $T:X\rightarrow Y$ un op\'{e}rateur born\'{e}. Alors $T^{\ast
\ast }$ est ponctuellement pr\'{e}faiblement absolument continu
\end{remark}

Preuve. 

En effet,

consid\'{e}rons $A\in \Sigma ,$ $f^{\ast \ast }\in X^{\ast \ast }$ et $%
y^{\ast }\in Y^{\ast }.$ Remarquons que $\left\langle y^{\ast },T^{\ast \ast
}(f^{\ast \ast }\mathcal{X}_{A})\right\rangle =\left\langle \Delta _{X^{\ast
}}^{\ast }(A,T^{\ast }(y^{\ast })),f^{\ast \ast }\right\rangle ,$ donc la
mesure $A\rightarrow \left\langle y^{\ast },T^{\ast \ast }(f^{\ast \ast }%
\mathcal{X}_{A}))\right\rangle $ est $\mu -$d\'{e}nombrablement additive,
car $X^{\ast }$ est absolument continu.$\blacksquare $

\begin{corollary}
\label{ri}Soient $X$ un espace de fonctions au sens large sur $(\Omega
,\Sigma ,\mu ,\Delta _{X})$. Alors $K(X,Y)$ est absolument continu si et
seulement si $X^{\ast }$ est absolument continu.
\end{corollary}

D\'{e}monstration.

Supposons que $K(X,Y)$ est absolument continu. Pour $f^{\ast }\in X^{\ast },$
on d\'{e}finit $U_{f^{\ast }}:X\rightarrow Y,$ par $U_{f^{\ast
}}(f)=\left\langle f,f^{\ast }\right\rangle y_{0},$ $f\in X,$ o\`{u} $%
\left\Vert y_{0}\right\Vert _{Y}=1.$ Fixons $f^{\ast }\in X^{\ast }$. Comme $%
U_{f^{\ast }}$ est absolument continu, pour tout $\varepsilon >0,$ il existe 
$\delta >0$ telle que sup$\left\{ \left\Vert U_{f^{\ast }}(f\mathcal{X}%
_{A})\right\Vert _{Y};\text{ }f\in B_{X}\right\} <\varepsilon ,$ si $\mu
(A)<\delta .$ Choisissons $y^{\ast }\in Y^{\ast }$ tel que $\left\vert
\left\langle y_{0},y^{\ast }\right\rangle \right\vert =1.$ Pour tout $f\in
B_{X}$, on a%
\begin{eqnarray*}
\left\Vert U_{f^{\ast }}(f\mathcal{X}_{A})\right\Vert _{Y} &\geq &\left\vert
\left\langle U_{f^{\ast }}(f\mathcal{X}_{A}),y_{{}}^{\ast }\right\rangle
\right\vert \\
&=&\left\vert \left\langle f\mathcal{X}_{A}),f_{{}}^{\ast }\right\rangle
\right\vert \left\vert \left\langle y_{0},y^{\ast }\right\rangle \right\vert
=\left\vert \left\langle f,\Delta _{X^{\ast }}^{\ast }(A,f^{\ast
})\right\rangle \right\vert ,
\end{eqnarray*}

par cons\'{e}quent $\left\Vert \Delta _{X^{\ast }}(A,f^{\ast })\right\Vert
<\varepsilon ,$ si $\mu (A)<\delta ;$

Inversement, supposons que $X^{\ast }$ est absolument continu. Soit $T\in
K(X,Y).$ D'apr\`{e}s la remarque \ref{vo}, $T$ est ponctuellement pr\'{e}%
faiblement absolument continu. En appliquant la proposition \ref{uv}, on
voit que $T$ est absolument continu.$\blacksquare $

\begin{remark}
\label{dn}Dans \cite{Ni}, on introduit les op\'{e}rateurs absolument
continus, pour cette d\'{e}finition tout op\'{e}rateur compact est
absolument continu, si $X$ ne contient pas $\ell ^{1}$ isomorphiquement.
\end{remark}

\begin{corollary}
\label{dm}Soient $X$ un espace de fonctions sur $(\Omega ,\Sigma ,\mu
,\Delta _{X})$, $Y$ un espace de Banach et $T:X\rightarrow Y$ un op\'{e}%
rateur compact et ponctuellement faiblement absolument continu. Si $Y$ a la
propri\'{e}t\'{e} de l'approximation born\'{e}e, alors $T$ est absolument
continu.
\end{corollary}

D\'{e}monstration.

Comme $Y$  a la propri\'{e}t\'{e} de l'approximation born\'{e}e, il existe
une suite g\'{e}n\'{e}ralis\'{e}e $(U_{i})_{i\in I}$ d'op\'{e}rateurs du
rang finis:$Y\rightarrow Y$, telle que $\sup_{i\in I}\left\Vert
U_{i}\right\Vert <+\infty $ et $U_{i}\rightarrow I_{Y}$ uniform\'{e}ment sur
tout compact de $Y.$ Pour tout $i\in I,$ notons $T_{i}=U_{i}\circ T.$ Comme $%
T_{i}\rightarrow T$ dans $\mathcal{L}(X,Y),$ il suffit de montrer que $T_{i}$
est absolument continu, d'apr\`{e}s la proposition \ref{uv}, il suffit de
montrer que $(T_{i})^{\ast \ast }$ est ponctuellement pr\'{e}faiblement
absolument continu$.$

Fixons $i\in I.$ L'op\'{e}rateur $T_{i}$ est du rang fini, il existe donc
une suite de vecteurs $(f_{k}^{\ast })_{k\leq n}$ dans $X^{\ast }$ et une
suite de vecteurs $(y_{k})_{k\leq n}$ dans $Y$ tels que $T_{i}(f)=\underset{%
k\leq n}{\dsum }\left\langle f,f_{k}^{\ast }\right\rangle y_{k},$ $f\in X.$

Consid\'{e}rons $(A_{k})_{k\geq 0}$ une suite de sous-ensembles mesurables
deux-\`{a}-deux disjoints et $y^{\ast }\in Y^{\ast }$. Comme $T_{i}$ est
ponctuellement faiblement absolument continu, d'apr\`{e}s la remarque \ref%
{dd}, $\underset{k\geq 0}{\dsum }\left\langle f,(T_{i}\mathcal{X}%
_{A_{k}})^{\ast }(y^{\ast })\right\rangle =\underset{k\geq 0}{\dsum }%
\left\langle T_{i}(f\mathcal{X}_{A_{k}}),y^{\ast }\right\rangle =$ $%
\left\langle T_{i}(f\mathcal{X}\underset{k\geq 0}{_{\cup }}A_{k}),y^{\ast
}\right\rangle =\left\langle f,(T_{i}\mathcal{X}_{\underset{k\geq 0}{\cup
A_{k}}})^{\ast }y^{\ast }\right\rangle .$ Il en r\'{e}sulte que la s\'{e}rie 
$\underset{k\geq 0}{\dsum }(T_{i}\mathcal{X}_{A_{k}})^{\ast }y^{\ast }$
converge pr\'{e}faiblement vers $(T_{i}\mathcal{X}\underset{k\geq 0}{_{\cup }%
}A_{k})^{\ast }y^{\ast }$ dans $X^{\ast }$. D'autre part, l'op\'{e}rateur $%
(T_{i}\mathcal{X}_{A})^{\ast }$ est \`{a} valeurs dans un sous-espace de
dimension finie $F_{0}$ de $X^{\ast }$ ind\'{e}pendant de $A$, car $T_{i}%
\mathcal{X}_{A}$ est \`{a} valeurs dans un sous espace de $Y$ ind\'{e}%
pendant de $A,$ par cons\'{e}quent la s\'{e}rie $\underset{k\geq 0}{\dsum }%
(T_{i}\mathcal{X}_{A_{k}})^{\ast }y^{\ast }$ converge fortement vers $(T_{i}%
\mathcal{X}\underset{k\geq 0}{_{\cup }}A_{k})^{\ast }y^{\ast }$ dans $%
X^{\ast }.$

Soit $f^{\ast \ast }\in X^{\ast \ast }$. D'apr\`{e}s le lemme \ref{kg}, $%
\underset{k\geq 0}{\dsum }\left\langle y^{\ast },T_{i}^{\ast \ast }(\mathcal{%
X}_{A_{k}}f^{\ast \ast })\right\rangle =\underset{k\geq 0}{\dsum }%
\left\langle (T_{i}\mathcal{X}_{A_{k}})^{\ast }y^{\ast },f^{\ast \ast
}\right\rangle $, donc $T_{i}^{\ast \ast }$ est ponctuellement pr\'{e}%
faiblement absolument continu.$\blacksquare $

\bigskip Soient $X$ un espace de fonctions au sens large sur $(\Omega
,\Sigma ,\mu ,\Delta _{X})$, $Y$ un espace de Banach et $T:X\rightarrow Y$
un op\'{e}rateur born\'{e}. Notons $E_{T}$=$\left\{ T\mathcal{X}_{A};\text{ }%
A\in \Sigma \right\} $ et $F_{T}$ le sous-espace ferm\'{e} engendr\'{e} par $%
E_{T}$ dans $\mathcal{L}(X,Y).$

\begin{proposition}
\label{de}Soient $X$ un espace de fonctions au sens large sur $(\Omega
,\Sigma ,\mu ,\Delta _{X})$, $\mu $ une mesure born\'{e}e, $Y$ un espace de
Banach et $T:X\rightarrow Y$ un op\'{e}rateur born\'{e} ponctuellement
faiblement absolument continu. Supposons que pour tout $\xi \in
(F_{T}^{{}})^{\ast },$ il existe une suite $(\xi _{n})_{n\geq 0}$ dans $%
X\otimes Y^{\ast }$ telle que $\xi _{n}\underset{n\rightarrow \infty }{%
\rightarrow }\xi $ pr\'{e}faiblement dans $(F_{T})^{\ast }$. Alors $T$ est
absolument continu.
\end{proposition}

D\'{e}monstration.

Soit $\xi \in (F_{T})^{\ast }.$ Il existe une suite $(\xi _{n})_{n\geq 0}$
dans $X\otimes Y^{\ast }$ telle $\xi _{n}\underset{n\rightarrow \infty }{%
\rightarrow }\xi $ pr\'{e}faiblement. Fixons $n\in \mathbb{N}$. $\ $On d\'{e}%
finit la mesure $\nu _{n}$ par $\nu _{n}(A)=\left\langle T\mathcal{X}%
_{A},\xi _{n}\right\rangle .$ D'apr\`{e}s l'hypoth\`{e}se, $\lim_{\mu
(A)\rightarrow 0}\nu _{n}(A)=0$ pout tout $n\in \mathbb{N}$. D'autre part, $%
\lim_{n\rightarrow \infty }\nu _{n}(A)=\left\langle T\mathcal{X}_{A},\xi
\right\rangle $ pour tout $A\in \Sigma $, d'apr\`{e}s\cite[Cor.6,Chap.1-5]%
{DU}, $\lim_{\mu (A)\rightarrow 0}\nu _{n}(A)=0,$ uniform\'{e}ment en $n\in 
\mathbb{N}.$ Il en r\'{e}sulte que $\lim_{\mu (A)\rightarrow 0}\left\langle T%
\mathcal{X}_{A},\xi \right\rangle =0.$ Soit $(A_{k})_{k\geq 0}$ une suite
dans $\Sigma $ deux-\`{a}-deux disjoint. D'apr\`{e}s la remarque \ref{dd},
pour tout $M\subset \mathbb{N}$ $\underset{k\in M}{\dsum }\left\langle T%
\mathcal{X}_{A_{k}},\xi \right\rangle =\left\langle T\mathcal{X}_{\underset{%
k\in M}{\cup }A_{k}},\xi \right\rangle ,$ pour tout $\xi \in (F_{T})^{\ast }.
$ Il en r\'{e}sulte que la s\'{e}rie $\underset{k\geq 0}{\dsum }T\mathcal{X}%
_{A_{k}}$ converge inconditionnellement vers $T\mathcal{X}_{\underset{k\geq 0%
}{\cup A_{k}}}$ dans $\mathcal{L}(X,Y)$ \cite[Corol.6,Chap.1-4]{DU}$%
.\blacksquare $

\begin{corollary}
\label{ry}Soient $X$ un espace de fonctions au sens large sur $(\Omega
,\Sigma ,\mu ,\Delta _{X})$, $\mu $ une mesure born\'{e}e, $Y$ un espace de
Banach et $T:X\rightarrow Y^{\ast }$ un op\'{e}rateur born\'{e}
ponctuellement faiblement absolument continu. Supposons que $X,Y$ sont s\'{e}%
parables et que $X\overset{\wedge }{\otimes }Y$ ne contient pas $\ell ^{1}$
isomorphiquement. Alors $T$ est absolument continu.
\end{corollary}

D\'{e}monstration.

D'apr\`{e}s \cite[Chap.VIII-2,Coroll.2]{DU}-\cite{Sch}, ($X\overset{\wedge }{%
\otimes }Y)^{\ast }=\mathcal{L}(X,Y^{\ast }).$ Soit $\xi \in \left[ \mathcal{%
L}(X,Y^{\ast })\right] ^{\ast }=(X\overset{\wedge }{\otimes }Y)^{\ast \ast
}. $ Comme $X\overset{\wedge }{\otimes }Y$ ne contient pas $\ell ^{1},$ d'apr%
\`{e}s \cite{Ros}, il existe une suite $(\xi _{n})_{n0}$ dans $X\overset{%
\wedge }{\otimes }Y$ telle $\xi _{n}\underset{n\rightarrow \infty }{%
\rightarrow }\xi $ pr\'{e}faiblement. Pour conclure, il suffit d'appliquer
la proposition \ref{de}.$\blacksquare $

\begin{proposition}
\label{jo}Soient $X$ un espace de fonctions au sens large sur $(\Omega
,\Sigma ,\mu ,\Delta _{X})$, $\mu $ une mesure born\'{e}e, $Y$ un espace de
Banach et $T:X\rightarrow Y$ un op\'{e}rateur absolument continu. Supposons
que $F_{T}$ ne contient pas $c_{0}$. Alors $F_{T}$ est un espace $WCG.$
\end{proposition}

D\'{e}monstration.

Comme $T$ est absolument continu, d'apr\`{e}s le lemme \ref{ez}, $\nu _{T}$
est $\mu -$d\'{e}nombrablement additive. D'autre part, $F_{T}$ ne contient
pas $c_{0}$ et $\nu _{T}$ est born\'{e}e, d'apr\`{e}s \cite[Chap.1-4,Th.2]%
{DU}, $\nu _{T}$ est fortement additive. En appliquant le r\'{e}sultat de 
\cite[Chap.I-5,Coroll.3]{DU}, on voit que $E_{T}$ est faiblement compact,
donc $F_{T}$ est un espace $WCG.\blacksquare $

\begin{proposition}
\label{copy}Soient $X$ un espace de fonctions au sens large sur $(\Omega
,\Sigma ,\mu ,\Delta _{X})$, $\mu $ une mesure born\'{e}e, $Y$ un espace de
Banach et $T:X\rightarrow Y$ un op\'{e}rateur born\'{e} . Supposons que $%
F_{T}$ est un espace $WCG$ $\ $et $T$ est ponctuellement faiblement
absolument continu. Alors $T$ est absolument continu$.$
\end{proposition}

\bigskip D\'{e}monstration.

Comme $F_{T}$ est un espace $WCG,$ il existe un compact faible $K$ dans $%
F_{T}$ tel que l'espace ferm\'{e} engendr\'{e} par $K$ est \'{e}gale \`{a} $%
F_{T}.$ Soit $J:B_{(F_{T})^{\ast }}\rightarrow C(K)$ l'injection canonique.
Consid\'{e}rons maintenant $\xi $ dans la boule unit\'{e} de $(F_{T})^{\ast
}.$ D'apr\`{e}s le th\'{e}or\`{e}me de Hahn-Banach, il existe $\eta $ dans
la boule unit\'{e} de $\left[ \mathcal{L}(X,Y^{\ast \ast })\right] ^{\ast }$
qui prolonge $\xi .$ Comme $(X\overset{\wedge }{\otimes }Y^{\ast })^{\ast }$
est l'espace $\mathcal{L}(X,Y^{\ast \ast })$ \cite[Chap.VIII-2,Coroll.2]{DU}-%
\cite{Sch}, il existe une suite g\'{e}n\'{e}ralis\'{e}s $(\eta _{i})_{i\in
I} $ dans la boule unit\'{e} de $X\overset{\wedge }{\otimes }Y^{\ast }$
telle que $\eta _{i}\rightarrow \eta $ $\sigma (\left[ \mathcal{L}(X,Y^{\ast
\ast })\right] ^{\ast },\mathcal{L}(X,Y^{\ast \ast })).$ D\'{e}signons pour
tout $i\in I$ par $\xi _{i}$ la restriction de $\eta _{i}$ \`{a} $F_{T}.$
Nous avons alors $\xi _{i}\rightarrow \xi $ $\sigma ((F_{T}^{{}})^{\ast
},F_{T}),$ $J(\xi _{i})\rightarrow J(\xi )$ pour la topologie de la
convergence simple dans $C(K).$ D'apr\`{e}s \cite{Groth}, il existe une
suite $(\xi _{i_{n}})_{n\geq 0}$ telle que $J(\xi _{i_{n}})\underset{%
n\rightarrow \infty }{\rightarrow }J(\xi )$ dans ($C(K),\tau _{p}),$ ceci
implique que $\xi _{i_{n}}\underset{n\rightarrow \infty }{\rightarrow }\xi ,$
$\sigma ((F_{T})^{\ast },F_{T})$, d'apr\`{e}s la proposition \ref{de}, $T$
est absolument continu.$\blacksquare $

\begin{remark}
\label{xs}Dans la proposition \ref{copy} on peur remplacer $F_{T}$ est un
espace $WCG$ par l'existence d'un espace $WCG$ de $\mathcal{L}(X,Y)$ qui
contient $F_{T}.$
\end{remark}

\begin{proposition}
\label{dy}Soient $X$ un espace de fonctions au sens large sur $(\Omega
,\Sigma ,\mu ,\Delta _{X})$, $Y$ un espace de Banach et $T:X\rightarrow Y$
une application lin\'{e}aire. Supposons que $T$ soit $p-$sommant pour un $%
p\in \left[ 1,\infty \text{ }\right[ $ et $X^{\ast }$ est absolument
continu. Alors $T$ est absolument continu.
\end{proposition}

D\'{e}monstration.

Il suffit de montrer que $\nu _{T}$ est $\mu -$d\'{e}nomrablement additive,
d'apr\`{e}s le lemme \ref{ez}.

Soit $(A_{k})_{k\geq 0}$ une suite dans $\Sigma $ deux-\`{a}-deux disjoints.
Comme $T$ est $p-$sommant, il existe une mesure positive $G$ sur la boule
unit\'{e} de $X^{\ast }$ telle que $\left\Vert T(f)\right\Vert \leq
C\dint\limits_{B_{X^{\ast }}}\left\vert \left\langle f,f^{\ast
}\right\rangle \right\vert $\bigskip $dG(f^{\ast }).$ Donc pour tout $m\geq n
$ 
\begin{eqnarray}
\left\Vert \overset{m}{\underset{k=n}{\dsum }}\nu _{T}(A_{k})(f)\right\Vert 
&\leq &C\dint\limits_{B_{X^{\ast }}}\left\vert \overset{m}{\underset{k=n}{%
\dsum }}\left\langle f\mathcal{X}_{A_{k}}),f^{\ast }\right\rangle
\right\vert \bigskip dG(f^{\ast })  \notag \\
&=&C\dint\limits_{B_{X^{\ast }}}\left\vert \left\langle f,\underset{k=n}{%
\overset{m}{\dsum }}\Delta _{X^{\ast }}^{\ast }(A_{k},f^{\ast
})\right\rangle \right\vert \bigskip dG(f^{\ast })  \label{zz} \\
&\leq &C\left\Vert f\right\Vert \dint\limits_{B_{X^{\ast }}}\left\Vert 
\underset{k=n}{\overset{m}{\dsum }}\Delta _{X^{\ast }}^{\ast }(A_{k},f^{\ast
})\bigskip \right\Vert _{X^{\ast }}dG(f^{\ast }).  \notag
\end{eqnarray}%
Comme $X^{\ast }$ est absolument continu, $\left\Vert \overset{m}{\underset{%
k=n}{\dsum }}\Delta _{X^{\ast }}^{\ast }(A_{k},f^{\ast })\bigskip
\right\Vert _{X^{\ast }}\underset{m,n\rightarrow \infty }{\rightarrow }0$
pour toute $f^{\ast }\in B_{X^{\ast }}.$  En appliquant le th\'{e}or\`{e}me
de convergence domin\'{e}e, nous d\'{e}duisons que $\dint\limits_{B_{X^{\ast
}}}\left\Vert \overset{m}{\underset{k=n}{\dsum }}\Delta _{X^{\ast }}^{\ast
}(A_{k},f^{\ast })\bigskip \right\Vert _{X^{\ast }}dG(f^{\ast }).\underset{%
m,n\rightarrow \infty }{\rightarrow }0$ (Observons que pour tout $m\geq
n\left\Vert \underset{k=n}{\overset{m}{\dsum }}\Delta _{X^{\ast }}^{\ast
}(A_{k},f^{\ast })\bigskip \right\Vert _{X^{\ast }}=\left\Vert \Delta
_{X^{\ast }}^{\ast }(\underset{k=n}{\overset{m}{\cup }}(A_{k},f^{\ast
})\right\Vert _{X^{\ast }}\leq 1).$ Il en r\'{e}sulte d'apr\`{e}s (\ref{zz})
que $\left\Vert \overset{m}{\underset{k=n}{\dsum }}\nu
_{T}(A_{k})\right\Vert \underset{m,n\rightarrow \infty }{\rightarrow }0,$
c'est-\`{a}-dire que la s\'{e}rie $\underset{k\geq 0}{\dsum }\nu _{T}(A_{k})$
converge en norme dans $\mathcal{L}(X,Y).\blacksquare $

Soit $X,Y$ deux espaces de Banach. D\'{e}signons par $\Pi (X,Y)$ l'espaces
des op\'{e}rateurs $T:X\rightarrow Y$ tel que $T$ est $p-$sommant pour un $%
p\in \left[ 1,+\infty \right[ $ et par $\mathcal{L}_{0}(X,Y)$ l'adh\'{e}%
rence de $\Pi (X,Y)$ dans $\mathcal{L}(X,Y).$

\begin{corollary}
\label{fd}Soient $X$ un espace de fonctions au sens large sur $(\Omega
,\Sigma ,\mu ,\Delta _{X})$ et $Y$ un espace de Banach. Supposons que $%
X^{\ast }$ soit absolument continu. Alors $\mathcal{L}_{0}(X,Y)$ est
absolument continu.
\end{corollary}

\section{\protect\bigskip Espaces de fonctions au sens large pour les
espaces d'interpolation}

Soient $\overline{B}=(B_{0},B_{1})$ un couple d'interpolation au sens de 
\cite[chap. II]{Ber-Lof}, et $\theta \in \left] 0,1\right[ .$

Soit $S=\left\{ z\in \mathbb{C};\mathrm{\ }0\leq \func{Re}(z)\leq 1\right\} $
et $S^{0}=\left\{ z\in \mathbb{C};0<\func{Re}(z)<1\right\} .$ On d\'{e}signe
par $\mathcal{F}\mathfrak{(}\overline{B})$ l'espace des fonctions $F$ \`{a}
valeurs dans $B_{0}+B_{1},$ continues born\'{e}es sur $S$, holomorphes sur $%
S^{0},$ telles que, pour $j\in \left\{ 0,1\right\} $, l'application $\tau
\rightarrow $ $F(j+i\tau )$ est continue \`{a} valeurs dans $B_{j}$ et $%
\left\Vert F(j+i\tau )\right\Vert _{B_{j}}\rightarrow _{\left\vert \tau
\right\vert \rightarrow \infty }0.$ On le munit de la norme%
\begin{equation*}
\left\Vert F\right\Vert _{\mathcal{F}\mathfrak{(}\overline{B)}\mathrm{\ }%
}=\max \{sup_{\tau \in \mathbb{R}}\left\Vert F(i\tau )\right\Vert
_{B_{0}},sup_{\tau \in \mathbb{R}}\left\Vert F(1+i\tau )\right\Vert
_{B_{1}}\}.
\end{equation*}

Notons $\mathcal{F}_{0}(B_{0},B_{1})=\left\{ \overset{n}{\underset{k=0}{%
\dsum }}F_{k}\otimes b_{k};\text{ }b_{k}\in B_{0}\cap B_{1},F_{k}\in 
\mathcal{F}(\mathbb{C)},\text{ }n\in \mathbb{N}\right\} .$

\noindent L'espace $(B_{0},B_{1})_{\theta }=B_{\theta }=\left\{ F(\theta );%
\mathrm{\ }F\in \mathcal{F}\mathfrak{(}\overline{B})\right\} $ est de Banach 
\cite[th.4.1.2]{Ber-Lof}, pour la norme d\'{e}finie par

\begin{equation*}
\left\Vert a\right\Vert _{B_{\theta }}=\inf \left\{ \left\Vert F\right\Vert
_{\mathcal{F}\mathfrak{(}\overline{B})};\mathrm{\ }F(\theta )=a\right\} .
\end{equation*}

\noindent Toute $F\in \mathcal{F}\mathfrak{(}\overline{B})$ est repr\'{e}sent%
\'{e}e \`{a} partir de ses valeurs au bord en utilisant la mesure
harmonique, de densit\'{e} $Q_{0}(z,i\tau )$ et $Q_{1}(z,1+i\tau )$, $z\in
S^{0}$, $\tau \in \mathbb{R},$ \cite[section 4.5]{Ber-Lof}:

\begin{equation}
F(z)=\dint\limits_{\mathbb{R}}F(i\tau )Q_{0}(z,\tau )d\tau +\dint\limits_{%
\mathbb{R}}F(1+i\tau )Q_{1}(z,\tau )d\tau .  \label{rt}
\end{equation}

On note $\mathcal{G(}\overline{B})$ l'espace des fonctions $g$ \`{a} valeurs
dans $B_{0}+B_{1}$ , continues sur $S,$ holomorphes \`{a} l'interieur de $S,$
telles que

$(C)\quad \sup_{z\in S}\frac{\left\Vert g(z)\right\Vert _{B_{0}+B_{1}}}{%
(1+\left\vert z\right\vert )}<\infty .$

\noindent\ $\ \ (C^{\prime })$ $g(j+i\tau )-g(j+i\tau ^{\prime })\in B_{j}$, 
$\forall \tau ,\tau ^{\prime }\in \mathbb{R}$,$\ $\ $j\in \left\{
0,1\right\} $ et la quantit\'{e} suivante est finie: 
\begin{equation*}
\left\Vert g\right\Vert _{Q\mathcal{G(}\overline{B}\mathcal{)}}=\max \left[ 
\begin{array}{c}
\sup_{\tau \neq \tau ^{\prime }\in R}(\left\Vert (g(i\tau )-g(i\tau ^{\prime
}))/\tau -\tau ^{\prime }\right\Vert _{B_{0}}), \\ 
\sup_{\tau \neq \tau ^{\prime }\in R}(\left\Vert (g(1+i\tau )-g(1+i\tau
^{\prime }))/\tau -\tau ^{\prime }\right\Vert _{B_{1}})%
\end{array}%
\right] .
\end{equation*}

\noindent Ceci d\'{e}finit bien une norme sur $Q\mathcal{G}(\overline{B})$
le quotient de $\mathcal{G}(\overline{B})$ (par les constante \`{a} valeurs
dans $B_{0}+B_{1}).$

\noindent L'espace $(B_{0},B_{1})^{\theta }=B^{\theta }=\left\{ \mathrm{\ }%
g^{\prime }(\theta );\mathrm{\ }g\in \mathcal{G(}\overline{B})\mathrm{\ }%
\right\} $ est de Banach \cite[th.4.1.4]{Ber-Lof} pour la norme d\'{e}finie
par

\begin{equation*}
\left\Vert a\right\Vert _{B^{\theta }}=\inf \left\{ \left\Vert g\right\Vert
_{Q\mathcal{G}(\overline{B})};\mathrm{\ }g^{\prime }(\theta )=a\right\} .%
\mathrm{\ }
\end{equation*}

Soit $p\in \left[ 1,+\infty \right[ $. L'espace d'interpolation $B_{\theta
,p}$ \ est d\'{e}fini par 
\begin{equation*}
B_{\theta ,p}=\left\{ a\in B_{0}+B_{1};\mathrm{\ }\left\Vert a\right\Vert
_{B_{\theta ,p}}=\left[ \int_{\mathbb{R}^{+}}(K(a,t)/t^{\theta })^{p}dt/t%
\right] ^{1/p}<+\infty \right\}
\end{equation*}

\noindent o\`{u}

\begin{equation*}
K(a,t)=\inf \left\{ \left\Vert a_{0}\right\Vert _{B_{0}}+t\left\Vert
a_{1}\right\Vert _{B_{1}};\mathrm{\ }a=a_{0}+a_{1},\mathrm{\ }a_{j}\in B_{j}%
\mathrm{,}\text{ }j\in \left\{ 0,1\right\} \right\} .
\end{equation*}

\noindent ($B_{\theta ,p}$, $\left\Vert .\right\Vert _{B_{\theta ,p}})$ est
un espace de Banach \cite[th 3.4.2]{Ber-Lof}$.$

\begin{definition}
\label{ze}Soit $(B_{0},B_{1})$ un couple d'interpolation. Supposons que $%
B_{j}$ est un espace de fonctions au sens large sur $(\Omega ,\Sigma ,\mu
,\Delta _{j}),$ $\ j\in \left\{ 0,1\right\} .$ On dit que $(B_{0},B_{1})$
est un couple d'interpolation compatible avec $(\Omega ,\Sigma ,\mu ,(\Delta
_{0},\Delta _{1}))$ si $\Delta _{0}(A,f)=\Delta _{1}(A,f),$ pour tout $A\in
\Sigma $ et tout $f\in B_{0}\cap B_{1}.$
\end{definition}

Il est facile de montrer la proposition suivante:

\begin{proposition}
\label{rp}Soient $(B_{0},B_{1})$ un couple d'interpolation compatible avec $%
(\Omega ,\Sigma ,\mu ,(\Delta _{0},\Delta _{1}))$ et $f_{j},g_{j}\in B_{j},$ 
$j\in \left\{ 0,1\right\} .$ Supposons que $f_{0}+f_{1}=g_{0}+g_{1}.$ Alors $%
\Delta _{0}(A,f_{0})+\Delta _{1}(A,f_{1})=\Delta _{0}(A,g_{0})+\Delta
_{1}(A,g_{1})$ pour tout $A\in \Sigma .$
\end{proposition}

On d\'{e}finit $(\Delta _{0}+\Delta _{1})(A,f)=\Delta _{0}(A,f_{0})+\Delta
_{1}(A,f_{1}),$ o\`{u} $A\in \Sigma $ et $f=f_{0}+f_{1}\in B_{0}+B_{1}.$

On remarque que $B_{0}+B_{1}$ est un espace de fonctions au sens large sur $%
(\Omega ,\Sigma ,\mu ,\Delta _{0}+\Delta _{1}).$

\begin{example}
\bigskip \label{cn}Soient $X,Y$ deux espaces de fonctions au sens large sur $%
(\Omega ,\Sigma ,\mu ,\Delta _{X})$,$(\Omega ,\Sigma ,,\mu ,\Delta _{Y})$
respectivement et $i:X\rightarrow Y$ une injection continue. Supposons que $%
i(\Delta _{X})=\Delta _{Y}.$ Alors $(X,Y)$ est un couple d'interpolation
compatible avec $(\Omega ,\Sigma ,\mu ,(\Delta _{X},\Delta _{Y})).$
\end{example}

\begin{proposition}
\label{fh}Soient $(B_{0},B_{1})$ un couple d'interpolation compatible avec $%
(\Omega ,\Sigma ,\mu ,(\Delta _{0},\Delta _{1}))$ et $\theta \in \left] 0,1%
\right[ .$ Alors $A_{\theta }$ est un espace de fonctions au sens large sur $%
(\Omega ,\Sigma ,\mu ,\Delta _{\theta }),$ o\`{u} $\Delta _{\theta
}(A,f)=(\Delta _{0}+\Delta _{1})(A,f),$ $(A$,$f)\in \Sigma \times B_{\theta
}.$
\end{proposition}

\bigskip D\'{e}monstration.

Soit $(A,f)\in \Sigma \times B_{\theta }.$ Montrons que $\Delta _{\theta
}(A,f)\in B_{\theta }$ et $\left\Vert \Delta _{\theta }(A,f)\right\Vert
_{B_{\theta }}\leq \left\Vert f\right\Vert _{B_{\theta }}.$ Il existe $F\in 
\mathcal{F}(B_{0},B_{1})$ telle que $F(\theta )=f.$ L'op\'{e}rateur $%
u\rightarrow (\Delta _{0}+\Delta _{1})(A,u)$ est born\'{e}, donc
l'application $F_{A}:$ $z\in S\rightarrow (\Delta _{0}+\Delta
_{1})(A,F(z))\in B_{0}+B_{1}$ est holomorphe sur $S^{0}$. Comme $(\Delta
_{0}+\Delta _{1})(A,(F(j+i.))=\Delta _{j}(A,F(j+i.)),$ $j\in \left\{
0,1\right\} ,$ $F_{A}\in \mathcal{F}(B_{0},B_{1}).$ Il est clair que

$\left\Vert F_{A}\right\Vert _{\mathcal{F}(B_{0},B_{1})}\leq \left\Vert
F\right\Vert _{\mathcal{F}(B_{0},B_{1})}.$ Donc $\left\Vert \Delta _{\theta
}(A,f)\right\Vert _{B_{\theta }}\leq \left\Vert f\right\Vert _{B_{\theta
}}.\blacksquare $

\begin{proposition}
\bigskip \label{kl}Soient $(B_{0},B_{1})$ un couple d'interpolation
compatible avec $(\Omega ,\Sigma ,\mu ,(\Delta _{0},\Delta _{1}))$, $\theta
\in \left] 0,1\right[ $ et $p\in \left] 1,+\infty \right[ .$ Alors $%
A_{\theta ,p}$ est un espace de fonctions au sens large sur $(\Omega ,\Sigma
,\mu ,\Delta _{\theta ,p}),$ o\`{u} $\Delta _{\theta ,p}(A,f)=(\Delta
_{0}+\Delta _{1})(A,f),$ $(A,f)\in \Sigma $ $\times B_{\theta ,p}.$
\end{proposition}

D\'{e}monstration.

Soient $(A,$ $f)\in \Sigma \times B_{\theta ,p}$ et $\varepsilon ,t>0.$ Il
existe $f_{0}\in B_{0}$ et $f_{1}\in B_{1}$ tel que $f=f_{0}+f_{1},$ et $%
K(f,t)+\varepsilon >\left\Vert f_{0}\right\Vert _{B_{0}}+t\left\Vert
f_{1}\right\Vert _{B_{1}}.$ D'autre part, $\Delta _{\theta ,p}(A,f)=\Delta
_{0}(A,f_{0})+\Delta _{1}(A,f_{1}),$ donc 
\begin{eqnarray*}
&&K\left[ \Delta _{\theta ,p}(A,f),t\right] \left\Vert \leq \Delta
_{0}(A,f_{0})\right\Vert _{B_{0}}+t\left\Vert \Delta
_{1}(A,f_{1})\right\Vert _{B_{1}} \\
&\leq &\left\Vert f_{0}\right\Vert _{B_{0}}+t\left\Vert f\right\Vert
_{B_{1}}\leq K(f,t)+\varepsilon .
\end{eqnarray*}%
Il en r\'{e}sulte que $\Delta _{\theta ,p}(A,f)\in B_{\theta }$ et $%
\left\Vert \Delta _{\theta ,p}(A,f)\right\Vert _{B_{\theta ,p}}\leq
\left\Vert f\right\Vert _{B_{\theta ,p}}$.$\blacksquare $

\begin{remark}
\bigskip \label{xw}Soient $X,Y$ deux espaces de fonctions au sens large sur $%
(\Omega ,\Sigma ,\mu ,\Delta _{X})$,$(\Omega ,\Sigma ,,\mu ,\Delta _{Y})$
respectivement, $i:X\rightarrow Y$ une injection continue, $\theta \in \left]
0,1\right[ $ et $p\in \left] 1,+\infty \right[ .$ Supposons que $i(\Delta
_{X})=\Delta _{Y}.$ D'apr\`{e}s la proposition \ref{fh}, $(X,Y)_{\theta }$
est un espace de fonctions au sens large sur ($\Omega ,\Sigma ,\mu ,\Delta
_{\theta }).$ Nous remarquons que pour tout $0<\theta <\beta <1,$ $i(\Delta
_{\theta })=\Delta _{\beta },$ o\`{u} $i:(X,Y)_{\theta }\rightarrow
(X,Y)_{\beta }$ l'injection canonique.
\end{remark}

\begin{proposition}
\bigskip \label{nb}Soient $X,Y$ deux espaces de fonctions au sens large sur $%
(\Omega ,\Sigma ,\mu ,\Delta _{X})$,$(\Omega ,\Sigma ,\mu ,\Delta _{Y})$
respectivement et $i:X\rightarrow Y$ une injection absolument continu,
d'image dense. Si $i(\Delta _{X})=\Delta _{Y}$, alors pour tout $0<\theta
<\beta <1,$ $i:(X,Y)_{\theta }\rightarrow (X,Y)_{\beta }$ est absolument
continu$.$
\end{proposition}

D\'{e}monstration.

\emph{Etape 1: }\ Soit $\beta \in \left] 0,1\right[ .$ Montrons que $%
i:(X,Y)_{\beta }\rightarrow Y$ est absolument continu.

Soit $\varepsilon >0.$ Comme $i:X\rightarrow Y$ est absolument continu, il
existe $\delta >0$ tel que si $\mu (A)<\delta ,$ alors $\left\Vert \nu
_{i}(A))\right\Vert _{\mathcal{L}(X,Y)}<\varepsilon ^{\frac{1}{1-\theta }}.$

Consid\'{e}rons $f\in X$ tel que $\left\Vert f\right\Vert _{(X,Y)_{\theta
}}<1.$ Il existe $F\in \mathcal{F}_{0}(X,Y)$ v\'{e}rifiant $F(\theta )=f$ et 
$\left\Vert F\right\Vert _{\mathcal{F}(X,Y)}<1.$ Fixons un ensemble
mesurable $A$ de $\Omega $ tel $\mu (A)<\delta .$ La fonction $z\in
S\rightarrow F(z)\mathcal{X}_{A}=(\Delta _{0}+\Delta _{1})(A,F(z))$ est dans 
$\mathcal{F}(X,Y)\subset \mathcal{F}(Y),$ donc d'apr\`{e}s \cite[Lemme 4.3.2]%
{Ber-Lof} 
\begin{eqnarray*}
&&\left\Vert i\left[ F(\theta )\mathcal{X}_{A}\right] \right\Vert
_{Y}=\left\Vert \Delta _{Y}(A,iF(\theta ))\right\Vert _{Y} \\
&\leq &\left[ \dint\limits_{\mathbb{R}}\left\Vert F(i\tau )\mathcal{X}%
_{A}\right\Vert _{Y}\frac{Q(\theta ,i\tau )}{1-\theta }d\tau \right]
^{1-\theta }\times \left[ \dint\limits_{\mathbb{R}}\left\Vert F(1+i\tau )%
\mathcal{X}_{A}\right\Vert _{Y}\frac{Q(\theta ,1+i\tau }{\theta })d\tau %
\right] ^{\theta } \\
&\leq &\left[ \dint\limits_{\mathbb{R}}\left\Vert \nu _{i}(A))\right\Vert _{%
\mathcal{L}(X,Y)}\times \left\Vert F(i\tau )\right\Vert _{X}\frac{Q(\theta
,i\tau )}{1-\theta }d\tau \right] ^{1-\theta }\times \\
&&\left[ \dint\limits_{\mathbb{R}}\left\Vert F(1+i\tau \right\Vert _{Y}\frac{%
Q(\theta ,1+i\tau )}{\theta }d\tau \right] ^{\theta }\leq \varepsilon .
\end{eqnarray*}

$\emph{Etape}$ $\emph{2:}$ Soit $\beta \in \left] 0,1\right[ .$ Montrons que 
$i:X\rightarrow (X,Y)_{\beta }$ est absolument continu.

D'apr\`{e}s le lemme \ref{rh} et la remarque \ref{nh}, $i^{\ast }:Y^{\ast
}\rightarrow X^{\ast }$ est absolument continu, donc d'apr\`{e}s l'\'{e}tape
1 $i^{\ast }:(X^{\ast },Y^{\ast })_{\beta }\rightarrow X^{\ast }$ est
absolument continu. En r\'{e}appliquant le lemme \ref{rh} et la remarque \ref%
{nh}, on voit que $i^{\ast \ast }:X^{\ast \ast }\rightarrow \left[ (X^{\ast
},Y^{\ast })_{\beta }\right] ^{\ast }$ est absolument continu. D'autre part,
d'apr\`{e}s \cite{Da}, la restriction de $i^{\ast \ast }$ \`{a} $X$ est \`{a}
valeurs dans $(X,Y)_{\beta }$ et $(X,Y)_{\beta }$ est un sous-espace ferm%
\'{e} de $\left[ (X^{\ast },Y^{\ast })_{\beta }\right] ^{\ast }$. Donc $%
i:X\rightarrow (X,Y)_{\beta }$ est absolument continu.

$\emph{Etape}$ $\emph{3:}$ Soit $\theta ,\beta \in \left] 0,1\right[ $ tel
que $\theta <\beta .$ Montrons que $i:(X,Y)_{\theta }\rightarrow
(X,Y)_{\beta }$ est absolument continu.

D'apr\`{e}s l'\'{e}tape 1, $i:(X,Y)_{\theta }\rightarrow $\ $Y$ est
absolument continu, en appliquant l'\'{e}tape 2, on voit que $%
i:(X,Y)_{\theta }\rightarrow \left[ (X,Y)_{\theta },Y\right] _{\eta }$

est absolument continu, pour tout $\eta \in \left] 0,1\right[ .$

Choisissons $(1-\eta )\theta =\beta .$ D'apr\`{e}s le th\'{e}or\`{e}me de r%
\'{e}it\'{e}ration \cite[Th.4.6.1]{Ber-Lof}, $\left[ (X,Y)_{\theta },Y\right]
_{\eta }=(X,Y)_{\beta }.$ On en d\'{e}duit que $i:(X,Y)_{\theta }\rightarrow
(X,Y)_{\beta }$ est ponctuellement absolument continu.$\blacksquare $

\begin{lemma}
\label{my}Supposons que $i$: $X\rightarrow Y$ soit une injection continue
d'image dense. Alors pour tout $0<\theta <1$ et tout $1<p<+\infty $ $%
(X,Y)_{\theta ,p}$ est un sous-espace ferm\'{e} de $\left[ (X^{\ast
},Y^{\ast })_{\theta ,p^{\prime }}\right] ^{\ast }$, ou $p^{\prime }$ est le
conjugu\'{e} de $p.$
\end{lemma}

D\'{e}monstration.

Il est \'{e}vident que l'injection :$((X,Y)_{\theta ,p}\rightarrow \left[
(X^{\ast },Y^{\ast })_{\theta ,p^{\prime }}\right] ^{\ast }$ est continue.
D'autre part, d'apr\`{e}s \cite[Th.3.7.1]{Ber-Lof}, il existe une constante $%
C>0,$ telle que $\left\Vert f^{\ast }\right\Vert _{(X^{\ast },Y^{\ast
})_{\theta ,p^{\prime }}}\leq C\left\Vert f^{\ast }\right\Vert _{\left[
(X,Y)_{\theta ,p}\right] ^{\ast }}$ pour tout $f^{\ast }\in \left[
(X,Y)_{\theta ,p}\right] ^{\ast }.$ Soient $f\in (X,Y)_{\theta ,p}$ et $%
\varepsilon >0.$ Il existe $f^{\ast }$ dans la boule unit\'{e} de $\left[
(X,Y)_{\theta ,p}\right] ^{\ast }$ tel que $\left\Vert f\right\Vert
_{(X,Y)_{\theta ,p}}\leq \left\vert \left\langle f,f^{\ast }\right\rangle
\right\vert +\varepsilon .$ D'apr\`{e}s ce qui pr\'{e}c\`{e}de, $f^{\ast
}\in (X^{\ast },Y^{\ast })_{\theta ,p^{\prime }}$ et $\left\Vert f^{\ast
}\right\Vert _{(X^{\ast },Y^{\ast })_{\theta ,p^{\prime }}}\leq C,$ d'o\`{u}
le lemme.$\blacksquare $

Par un argument analogue \`{a} celui de la proposition \ref{nb} (en
utilisant le lemme \ref{my}) on montre:

\begin{proposition}
\label{vv}Soient $X,Y$ deux espaces de fonctions au sens large sur $(\Omega
,\Sigma ,\mu ,\Delta _{X})$,$(\Omega ,\Sigma ,\mu ,\Delta _{Y})$
respectivement et $i:X\rightarrow Y$ une injection absolument continu
d'image dense. Supposons que $i(\Delta _{X})=\Delta _{Y}$. Alors pour tout $%
0<\theta <\beta <1$ et tout $1<p<+\infty ,$ $i:(X,Y)_{\theta ,p}\rightarrow
(X,Y)_{\beta ,p}$ est absolument continu$.$
\end{proposition}

\begin{lemma}
\label{mo}Soient $X$ un espace de fonctions au sens large sur $(\Omega
,\Sigma ,\mu ,\Delta _{X}),$ $Y$ un espace de Banach et $i:X\rightarrow Y$
une injection ponctuellement faiblement absolument continu, $\theta \in %
\left] 0,1\right[ $ et $p\in \left] 1,+\infty \right[ $. Alors $%
i:X\rightarrow (X,Y)_{\theta ,p}$ est ponctuellement faiblement absolument
continu.
\end{lemma}

D\'{e}monstration.

Soient $f\in X$, $x^{\ast }\in \left[ (X,Y)_{\theta ,p}\right] ^{\ast }$ et $%
\varepsilon ,\varepsilon ^{\prime }>0$ tel que $\left\Vert i\right\Vert
\left\Vert f\right\Vert \varepsilon ^{\prime }+\varepsilon ^{\prime }\leq
\varepsilon $ $.$ D'apr\`{e}s \cite[Th.3.7.1]{Ber-Lof} $\left[ (X,Y)_{\theta
,p}\right] ^{\ast }=(X^{\ast },Y^{\ast })_{\theta ,p^{\prime }}$. D'autre
part, le th\'{e}or\`{e}me 3.4.2 de \cite{Ber-Lof} nous montre que $Y^{\ast }$
est dense dans $\left[ (X,Y)_{\theta ,p}\right] ^{\ast },$ par cons\'{e}%
quent il existe $z^{\ast }\in Y^{\ast }$ tel que $\left\Vert x^{\ast
}-z^{\ast }\right\Vert $ $_{\left[ (X,Y)_{\theta ,p}\right] ^{\ast
}}<\varepsilon ^{\prime }.$ D'apr\`{e}s l'hypoth\`{e}se, il existe $\delta
>0 $ tel que si $\mu (A)<\delta ,$ $\left\vert \left\langle i(f\mathcal{X}%
_{A}),z^{\ast })\right\rangle \right\vert <\varepsilon ^{\prime }.$ Soit $%
A\in \Sigma $ tel que $\mu (A)<\delta .$ On a alors%
\begin{eqnarray*}
\left\vert \left\langle i(f\mathcal{X}_{A}),x^{\ast }\right\rangle
\right\vert &=&\left\vert \left\langle i(f\mathcal{X}_{A}),x^{\ast }-z^{\ast
}+z^{\ast }\right\rangle \right\vert \leq \\
\left\vert \left\langle i(f\mathcal{X}_{A}),x^{\ast }-z^{\ast }\right\rangle
\right\vert +\left\vert \left\langle i(f\mathcal{X}_{A}),z^{\ast
}\right\rangle \right\vert &\leq & \\
\left\Vert i(f\mathcal{X}_{A})\right\Vert _{(X,Y)_{\theta ,p}}\times
\left\Vert x^{\ast }-z^{\ast }\right\Vert _{\left[ (X,Y)_{\theta ,p}\right]
^{\ast }}+\varepsilon ^{\prime } &\leq & \\
\left\Vert i\right\Vert \left\Vert f\right\Vert \varepsilon ^{\prime
}+\varepsilon ^{\prime } &\leq &\varepsilon .
\end{eqnarray*}

Il en r\'{e}sulte que $i:X\rightarrow (X,Y)_{\theta ,p}$ est faiblement
absolument continu.$\blacksquare $

\begin{lemma}
\label{io}Soient $X,Y$ deux espaces de fonctions au sens large sur $(\Omega
,\Sigma ,\mu ,\Delta _{X}),(\Omega ,\Sigma ,\mu ,\Delta _{Y})$ respectivement%
$,$ $\theta \in \left] 0,1\right[ ,$ $i:X\rightarrow Y$ une injection
ponctuellement faiblement absolument continu et $p\in \left] 1,+\infty %
\right[ .$ Supposons que $i(\Delta _{X})=\Delta _{Y}$. Alors $%
i:(X,Y)_{\theta ,p}\rightarrow Y$ est ponctuellement faiblement absolument
continu.
\end{lemma}

D\'{e}monstration.

Soient $f\in (X,Y)_{\theta ,p},$ $y^{\ast }\in Y^{\ast }$ et $\varepsilon
,\varepsilon ^{\prime }>0$ tel que $\varepsilon ^{\prime }+\left\Vert
y^{\ast }\right\Vert \left\Vert i\right\Vert \varepsilon ^{\prime }\leq
\varepsilon .$ Il existe $f_{1}\in X$ tel que $\left\Vert f-f_{1}\right\Vert
_{(X,Y)_{\theta ,p}}<\varepsilon ^{\prime }.$ D'apr\`{e}s l'hypoth\`{e}se il
existe $\delta >0,$ tel que si $\mu (A)<\delta $ $\left\vert \left\langle
i(f_{1}\mathcal{X}_{A}),y^{\ast }\right\rangle \right\vert <\varepsilon
^{\prime },.$ Choisissons $A\in \Sigma $ tel que $\mu (A)<\delta .$ On a
alors $\left\vert \left\langle i(f\mathcal{X}_{A}),y^{\ast }\right\rangle
\right\vert \leq \left\vert \left\langle i(f_{1}\mathcal{X}_{A}),y^{\ast
}\right\rangle \right\vert +\left\Vert y^{\ast }\right\Vert \left\Vert
i(f-f_{1})_{A}\right\Vert _{Y}<\varepsilon ^{\prime }+\left\Vert
i\right\Vert \left\Vert y^{\ast }\right\Vert \varepsilon ^{\prime }\leq
\varepsilon .\blacksquare $

Par un argument analogue \`{a} celui de la proposition \ref{nb} (en
utilisant le th\'{e}or\`{e}me de r\'{e}it\'{e}ration \cite[Th.3.5.3]{Ber-Lof}
et les lemmes \ref{mo}, \ref{io}) on tire le corollaire suivant:

\begin{corollary}
\label{tp}Soient $X,Y$ deux espaces de fonctions au sens large sur $(\Omega
,\Sigma ,\mu ,\Delta _{X}),(\Omega ,\Sigma ,\mu ,\Delta _{Y})$
respectivement, $0<\theta <\beta <1$, $\ p\in \left] 1,+\infty \right[ $ et $%
i:X\rightarrow Y$ une injection ponctuellement faiblement absolument
continu. Supposons que $i(\Delta _{X})=\Delta _{Y}$. Alors $i:(X,Y)_{\theta
,p}\rightarrow (X,Y)_{\beta ,p}$ est ponctuellement faiblement absolument
continu.
\end{corollary}

\begin{proposition}
\bigskip \label{jk}Soient $(B_{0},B_{1}),(C_{0},C_{1})$ deux couples
d'interpolation compatibles avec $(\Omega ,\Sigma ,,\mu ,(\Delta _{0},\Delta
_{1})),(\Omega ,\Sigma ,\mu ,(\Delta _{2,}\Delta _{3}))$ respectivement, $%
\theta \in \left] 0,1\right[ $ \ et $T:B_{j}\rightarrow C_{j}$ un op\'{e}%
rateur born\'{e}, $j\in \left\{ 0,1\right\} $ (avec $T_{0_{\mid B_{0}\cap
B_{1}}}=T_{1_{\mid B_{0}\cap B_{1}}})$. Supposons que $T:B_{0}\rightarrow
C_{0}$ est absolument continu. Alors $T:B_{\theta }\rightarrow C_{\theta }$
est absolument continu.
\end{proposition}

D\'{e}monstration.

Soit $A\in \Sigma .$ Observons que $T\mathcal{X}_{A}:B_{j}\rightarrow C_{j}$
est born\'{e}, d'apr\`{e}s \cite[Th.4.1.2]{Ber-Lof} $\left\Vert T\mathcal{X}%
_{A}\right\Vert _{B_{\theta }\rightarrow C_{\theta }}\leq \left[ \left\Vert T%
\mathcal{X}_{A}\right\Vert _{B_{0}\rightarrow C_{0}}\right] ^{1-\theta
}\times \left[ \left\Vert T\mathcal{X}_{A}\right\Vert _{B_{1}\rightarrow
C_{1}}\right] ^{\theta }.$ Donc $T:B_{\theta }\rightarrow C_{\theta }$ est
absolument continu.$\blacksquare $

\begin{proposition}
\label{vh}Soient $(B_{0},B_{1})$ un couple d'interpolation compatible avec $%
(\Omega ,\Sigma ,\mu ,(\Delta _{0},\Delta _{1}))$, $\theta \in \left] 0,1%
\right[ $. Supposons que $B_{0}$ est absolument continu. Alors $B_{\theta }$
est absolument continu.
\end{proposition}

D\'{e}monstration.

D'apr\`{e}s la remarque \ref{lr}, il suffit de montrer que $B_{0}\cap B_{1}$
est absolument continue dans $B_{\theta }.$ Pour cela, soient $f\in
B_{0}\cap B_{1}$ et $F\in \mathcal{F}_{0}(B_{0},B_{1})$ tels que $F(\theta
)=f.$ Posons $F_{A}(z)=(\Delta _{0}+\Delta _{1})(A,F(z)),$ $A\in \Sigma ,$ $%
z\in S,.$ Il est clair que $F_{A}\in \mathcal{F}(B_{0},B_{1}),$ d'apr\`{e}s 
\cite[Lemme.4.3.2]{Ber-Lof}, pour tout $A\in \Sigma $ nous avons 
\begin{eqnarray}
\left\Vert F_{A}(\theta )\right\Vert _{B_{\theta }} &\leq &\left[
\dint\limits_{\mathbb{R}}\left\Vert \Delta _{0}(A,F(i\tau )\right\Vert
_{B_{0}}\frac{Q(\theta ,i\tau )}{1-\theta }d\tau \right] ^{1-\theta }  \notag
\\
&&\times \left[ \dint\limits_{\mathbb{R}}\left\Vert \Delta _{1}(A,F(1+i\tau
)\right\Vert _{B_{1}}\frac{Q(\theta ,1+i\tau )}{\theta }d\tau \right]
^{\theta }.  \label{d}
\end{eqnarray}

Soit $(A_{n})_{n\geq 0}$ une suite dans $\Sigma $ telle que $\mu (A_{n})%
\underset{n\rightarrow \infty }{\rightarrow }0.$ Comme $B_{0}$ est
absolument continu pour tout $\tau \in \mathbb{R}$ $\Delta
_{0}(A_{n},F(i\tau ))\underset{n\rightarrow \infty }{\rightarrow }0,$ en
appliquant le th\'{e}or\`{e}me de convergence domin\'{e}e, on voit que$\left[
\dint\limits_{\mathbb{R}}\left\Vert \Delta _{0}(A_{n},F(i\tau )\right\Vert
_{B_{0}}\frac{Q(\theta ,i\tau )}{1-\theta }d\tau \right] \underset{%
n\rightarrow \infty }{\rightarrow }0.$ Il en r\'{e}sulte d'apr\`{e}s (\ref{d}%
) que $\Delta _{\theta }(A_{n},f)=F_{A_{n}}(\theta )\underset{n\rightarrow
\infty }{\rightarrow }0.\blacksquare $

\begin{definition}
\label{bm}Soient $X,Y$ deux espaces de fonctions sur $(\Omega ,\Sigma ,\mu
,\Delta _{X}),$ $(\Omega ,\Sigma ,\mu ,\Delta _{Y})$ respectivement. On dit
que $X,Y$ sont isom\'{e}triques au sens large, s'il existe un op\'{e}rateur
isom\'{e}trie surjectif $T:X\rightarrow Y$ v\'{e}rifiant $T(\Delta
_{X})=\Delta _{Y}.$
\end{definition}

\begin{remark}
\label{kh}Soient $(B_{0},B_{1})$ un couple d'interpolation compatible avec $%
(\Omega ,\Sigma ,\mu ,(\Delta _{0},\Delta _{1}))$, $\theta ,\alpha ,\beta
,\eta \in \left] 0,1\right[ $ tels que $\theta =(1-\eta )\alpha +\eta \beta
. $ Alors $A_{\theta }=(A_{\alpha },A_{\beta })_{\eta }$ isom\'{e}triquement
au sens large.
\end{remark}

Preuve. En effet,

Il est facile de voir que $(A_{\alpha },A_{\beta })$ est compatible avec $%
(\Omega ,\Sigma ,\mu ,(\Delta _{\alpha },\Delta _{\beta }))$. Consid\'{e}%
rons $\widetilde{\Delta }_{\eta }$ l'application qui d\'{e}finie l'espace de
fonctions $(B_{\alpha },B_{\beta })_{\eta }$ par rapport \`{a} $(\Omega
,\Sigma ,\mu ,(\Delta _{\alpha },\Delta _{\beta }))$, $f\in B_{0}\cap B_{1}$
et $A\in \Sigma $. Il est clair que $\widetilde{\Delta }_{\eta }(A,f)=\Delta
_{\alpha }(A,f)=\Delta _{\beta }(A,f)=\Delta _{0}(A,f)=\Delta _{\theta
}(A,f).$ D'autre part, d'apr\`{e}s le th\'{e}or\`{e}me de r\'{e}it\'{e}%
ration $B_{\theta }$ et $(B_{\alpha },B_{\beta })_{\eta }$ sont isom\'{e}%
triques et $\widetilde{\Delta }_{\eta }(A,f)=\Delta _{\theta }(A,f)$ pour
tout $A\in \Sigma $ et tout $f\in B_{\theta },$ car $B_{0}\cap B_{1}$ est
dense dans $(B_{\alpha },B_{\beta })_{\eta }=A_{\theta }$.$\blacksquare $

Pour tout $g\in \mathcal{G}(B_{0},B_{1})$ et tout $A\in \Sigma $ notons $%
g_{A}(z)=(\Delta _{0}+\Delta _{1})(A,g(z))$, $z\in S.$

\begin{lemma}
\label{cx}Soient $(B_{0},B_{1})$ un couple d'interpolation compatible avec $%
(\Omega ,\Sigma ,\mu ,(\Delta _{0},\Delta _{1}))$, $\theta \in \left] 0,1%
\right[ ,$ $A\in \Sigma $ et $g\in \mathcal{G}(B_{0},B_{1})$. Alors $%
g_{A}\in \mathcal{G}(B_{0},B_{1})$ et $\left\Vert g_{A}\right\Vert _{%
\mathcal{G}}\leq \left\Vert g\right\Vert _{\mathcal{G}}.$
\end{lemma}

D\'{e}monstration.

Il est clair qu $g_{A}$ est holomorphe sur \ $S^{0}$ et continue sur $S$ 
\`{a} valeurs dans $B_{0}+B_{1}$ $,$ car $(\Delta _{0}+\Delta _{1})(A,.)$
est un op\'{e}rateur born\'{e} sur $B_{0}+B_{1}.$

D'autre part, pour tout $\tau ,\tau ^{\prime }\in \mathbb{R}$ on a%
\begin{eqnarray}
g_{A}(j+i\tau )-g_{A}(j+i\tau ^{\prime }) &=&\Delta _{j}(A,g(j+i\tau
))-\Delta _{j}(A,g(j+i\tau ^{\prime })  \label{f} \\
&=&\Delta _{j}(A,g(j+i\tau ))-g(j+i\tau ^{\prime }))\in B_{j},\text{ \ }j\in
\left\{ 0,1\right\} .  \notag
\end{eqnarray}

La relation (\ref{f}) montre que $g_{A}$ v\'{e}rifie les condtions $%
C,C^{\prime }$ et $\left\Vert g_{A}\right\Vert _{\mathcal{G}}\leq \left\Vert
g\right\Vert _{\mathcal{G}}.\blacksquare $

\begin{remark}
\label{gk}Pour tout $\theta \in \left] 0,1\right[ $ et tout $g\in \mathcal{G}%
(B_{0},B_{1})$ on a $g_{A}(A,.)^{\prime }(\theta )=(\Delta _{0}+\Delta
_{1})(A,g^{\prime }(\theta )).$
\end{remark}

\bigskip D'apr\`{e}s le lemme \ref{cx} et la remarque \ref{gk}, on a la
proposition suivante:

\begin{proposition}
\label{ty}Soient $(B_{0},B_{1})$ un couple d'interpolation compatible avec $%
(\Omega ,\Sigma ,\mu ,(\Delta _{0},\Delta _{1}))$ et $\theta \in \left] 0,1%
\right[ .$ Alors $B^{\theta }$ est un espace de fonction au sens large sur $%
(\Omega ,\Sigma ,\mu ,\Delta ^{\theta }),$ o\`{u} $\Delta ^{\theta
}(A,f)=(\Delta _{0}+\Delta _{1})(A,f),$ $(A,f)\in \Sigma \times B^{\theta }.$
\end{proposition}

\begin{lemma}
\label{sq}Soient $(B_{0},B_{1})$ un couple d'interpolation compatible avec $%
(\Omega ,\Sigma ,\mu ,(\Delta _{0},\Delta _{1}))$ et $\theta \in \left] 0,1%
\right[ .$ Si $B_{0}\cap B_{1}$ est dense dans $B_{0}$ et $B_{1},$ alors $%
(B_{0}^{\ast },B_{1}^{\ast })$ est compactible avec $(\Omega ,\Sigma ,\mu
,(\Delta _{0}^{\ast },\Delta _{1}^{\ast }))$.
\end{lemma}

D\'{e}monstration.

Soient $f^{\ast }\in B_{0}^{\ast }\cap B_{1}^{\ast }$ et $A\in \Sigma .$
Remarquons que $\left\langle f,\Delta _{0}^{\ast }(A,f^{\ast })\right\rangle
=\left\langle f,\Delta _{1}^{\ast }(A,f^{\ast })\right\rangle $ pour tout $%
f\in B_{0}\cap B_{1}.$ D'autre part, d'apr\`{e}s \cite[Th.2.7.1]{Ber-Lof} $%
(B_{0}^{{}}\cap B_{1}^{{}})^{\ast }=B_{0}^{\ast }+B_{1}^{\ast },$ donc $%
\Delta _{0}^{\ast }(A,f^{\ast })=\Delta _{1}^{\ast }(A,f^{\ast })$ dans $%
B_{0}^{\ast }+B_{1}^{\ast }.\blacksquare $

\begin{proposition}
\label{wl}Soient $(B_{0},B_{1})$ un couple d'interpolation compatible avec $%
(\Omega ,\Sigma ,\mu ,(\Delta _{0},\Delta _{1}))$ et $\theta \in \left] 0,1%
\right[ $ Supposons que $B_{0}\cap B_{1}$ soit dense dans $B_{0}$ et $B_{1}.$
Alors $B_{\theta }^{\ast }=(B_{0}^{\ast },B_{1}^{\ast })^{\theta }$ isom\'{e}%
triquement au sens large.
\end{proposition}

D\'{e}monstration.

D'apr\`{e}s le th\'{e}or\`{e}me de dualit\'{e} \cite[Th.4.5.1]{Ber-Lof}, on
a isom\'{e}triquement $B_{\theta }^{\ast }=(B_{0}^{\ast },B_{1}^{\ast
})^{\theta }$. Consid\'{e}rons $f^{\ast }\in B_{\theta }^{\ast }$, $f\in
B_{0}\cap B_{1}$ et $\Delta ^{\theta }$ l'application qui d\'{e}finie
l'espace de fonctions ($B_{0}^{\ast },B_{1}^{\ast })^{\theta }.$ Il existe $%
g\in \mathcal{G}(B_{0}^{\ast },B_{1}^{\ast })$ tel que $g^{\prime }(\theta
)=f^{\ast }.$ Pour tout $A\in \Sigma $ on a 
\begin{equation*}
\left\langle f,\Delta _{\theta }^{\ast }(A,f^{\ast })\right\rangle
=\left\langle \Delta _{\theta }(A,f),g^{\prime }(\theta )\right\rangle .
\end{equation*}

D'autre part, $\Delta _{{}}^{\theta }(A,g^{\prime }(\theta ))=(\Delta
_{0}^{\ast }+\Delta _{1}^{\ast })(A,g^{\prime }(\theta ))$ et il existe $%
b_{j}^{\ast }\in B_{j}^{\ast }$ telle que $g^{\prime }(\theta )=b_{0}^{\ast
}+b_{1}^{\ast }.$ Ceci implique que%
\begin{eqnarray*}
\left\langle f,\Delta _{{}}^{\theta }(A,g^{\prime }(\theta )\right\rangle
&=&\left\langle f,(\Delta _{0}^{\ast }+\Delta _{1}^{\ast })(A,g^{\prime
}(\theta ))\right\rangle \\
&=&\left\langle f,\Delta _{0}^{\ast }(A,b_{0}^{\ast })+\Delta _{1}^{\ast
}(A,b_{1}^{\ast })\right\rangle =\left\langle \Delta _{0}(A,f),b_{0}^{\ast
}\right\rangle +\left\langle \Delta _{1}(A,f),b_{1}^{\ast }\right\rangle \\
&=&\left\langle \Delta _{\theta }(A,f),(b_{0}^{\ast }+b_{1}^{\ast
}\right\rangle =\left\langle \Delta _{\theta }(A,f),(g^{\prime }(\theta
)\right\rangle = \\
&=&\left\langle f,\Delta _{\theta }^{\ast }(A,g^{\prime }(\theta
))\right\rangle
\end{eqnarray*}

Il en r\'{e}sulte que $\left\langle f,\Delta _{\theta }^{\ast }(A,f^{\ast
}\right\rangle =\left\langle f,\Delta _{{}}^{\theta }(A,g^{\prime }(\theta
)\right\rangle $ pour $f\in B_{0}\cap B_{1}$ et $A\in \Sigma ,$ donc $\Delta
^{\theta }$ =$\Delta _{\theta }^{\ast }$ dans $B_{0}^{\ast }+B_{1}^{\ast
}=(B_{0}\cap B_{1})^{\ast },$ c'est-\`{a}-dire que $B_{\theta }^{\ast
}=(B_{0}^{\ast },B_{1}^{\ast })^{\theta }$ isom\'{e}triquement au sens large.%
$\blacksquare $

\bigskip


\begin{thebibliography}{Ber-Lof}
\bibitem[Ber]{Ber} \textbf{J. Bergh, }\textit{On the relation between the
two complex methods of interpolation, Indiana Univ. Math. J. 28, 775-777,
(1979).}

\bibitem[Ben-Sh]{Ben-Sh} \textbf{C. Bennet, R. Sharpley, }\textit{%
Interpolation of operators, Academie Press, (1988).}

\bibitem[Ber-Lof]{Ber-Lof} \textbf{J. Bergh, J. L\"{o}fstr\"{o}m, }\ \textit{%
Interpolation spaces an introduction, Springer-Verlag-Berlin Heidelberg New
York, (1976).}

\bibitem[Da]{Da} \textbf{M. Daher, }\textit{Une remarque sur les espaces
d'interpolation faiblement localement uniform\'{e}ment convexes},
arXiv:1206.4848.

\bibitem[Cal]{Cal} \textbf{A. P. Calder\'{o}n, }\textit{Intermediate spaces
and interpolation, the complex method, Studia Math. 24, 113-190, (1964).}

\bibitem[DU]{DU} \textbf{J. Diestel, J. J. Uhl}, \textit{Vector measures,
Math. Surveys 15 \ A.M.S, (1977).}

\bibitem[Groth]{Groth} \textbf{A. Grothendieck, }\textit{Crit\`{e}res de
compacit\'{e} dans les espaces fonctionnels g\'{e}n\'{e}raux, Amer. J. Math.
168-186, (1952).)}

\bibitem[Kalt]{Kalt} \textbf{N. J. Kalton,}\textit{\ Spaces of compact
operators, Math. Ann. 208, 267-278, (1974).}

\bibitem[Ni]{Ni} \textbf{Niculescu P. Consantin, }\textit{\ Absolute in
Banach spaces theory, Rev. Roum. Pure Appl. Vol. 24, 413-422, (1979). }

\bibitem[Ros]{Ros} \textbf{H. P. Rosenthal, }\textit{A characterzation of
spaces containing }$\ell ^{1}$\textit{, Proc. Nat. Sci. (U. S. A), Vol. 71,
No. 2, 411-2413, (1974).}

\bibitem[Sch]{Sch} \textbf{R. Schatten, }\textit{A theory of cross spaces,
Princetton Univ. Press, (1950).}
\end{thebibliography}
\end{document}